\documentclass[final]{siamltex}
\usepackage{amsfonts} 
\usepackage{a4wide}
\usepackage{epsfig}
\usepackage{graphicx}
\usepackage{subfig}  
\newcommand \rjpundemi {r_{j + 1/2}}
\newcommand \rjmundemi {r_{j - 1/2}}
\newcommand \rjjj {r_j} 
\newcommand \dVt {dV_{\widetilde g}}
\newcommand \DeltaVt {\Delta \widetilde V} 
\newcommand \tT         {\widetilde T} 
\newcommand \del            \partial 
\newcommand \RR             {\mathbb  R} 
\newcommand \St {\widetilde S}
\newcommand \la \langle
\newcommand \ra \rangle
\newcommand \xb {\overline x}
\newcommand \vb {\overline v} 
\newcommand{\auth}{\textsc}
\newcommand \be         {\begin{equation}}
\newcommand \ee         {\end{equation}}
\newcommand \R     {\mathbb{R}}
\newcommand \eps    {\epsilon}
\newcommand \lam     {\lambda}
\newcommand \om       {\omega}
\newcommand \omb     {\overline{\omega}}
\newcommand \Omb       {\overline{\Omega}} 
\newcommand{\dK}{\partial K}
\newcommand{\ekp}{e_{K}^{+}}
\newcommand{\ekm}{e_{K}^{-}} 
\newcommand \TT     {\mathcal{T}}
\newcommand \Tcal   {\TT} 
 
\title{Relativistic Burgers equations on curved spacetimes. Derivation and finite volume approximation} 

\author{Philippe G. L{\scriptsize e}Floch\thanks{Laboratoire Jacques--Louis Lions \& Centre National de la Recherche Scientifique, Universit\'e Pierre et Marie Curie (Paris 6), 4 Place Jussieu, 75252 Paris, France. {\tt (contact@philippelefloch.org, makhlof@ann.jussieu.fr)}}
\and Hasan Makhlof$^* \hskip-.15cm$  \and 
Baver Okutmustur\thanks{Department of Mathematics, Middle East Technical University, 06800 Ankara, Turkey. {\tt (baver@metu.edu.tr})}
\newline   
2000\textit{\ AMS Subject Class.} 35L65, 76N10, 83A05. 
}

\begin{document}

\maketitle

\begin{abstract} 
Within the class of nonlinear hyperbolic balance laws posed on a curved spacetime (endowed with a volume form), we identify a hyperbolic balance law that enjoys the same Lorentz invariance property as the one satisfied by the Euler equations of relativistic compressible fluids. This model is unique up to normalization and converges to the standard inviscid Burgers equation in the limit of infinite light speed. 
Furthermore, from the Euler system of relativistic compressible flows on a curved background, we derive both the standard inviscid Burgers equation and our relativistic generalizations. The proposed models are referred to as relativistic Burgers equations on curved spacetimes
and provide us with simple models on which numerical methods can be developed and analyzed.  
Next, we introduce a finite volume scheme for the approximation of discontinuous solutions to these relativistic Burgers equations. Our scheme is formulated geometrically and is consistent with the natural divergence form of the balance laws under consideration. It applies to weak solutions containing shock waves and, most importantly, is well balanced in the sense that it preserves steady solutions. Numerical experiments are presented which demonstrate the convergence of the proposed finite volume scheme and its relevance for computing entropy solutions on a curved background.
\end{abstract}


\section{Introduction}

In computational fluid dynamics, the (inviscid) Burgers equation 
$$
\del_t u + \del_x \big( u^2 /2\big) = 0, \qquad u=u(t,x), \quad t >0, \, x \in \RR,
$$
is an important toy model, which can be formally derived from the Euler equations of compressible fluids
$$
\del_t \rho + \del_x (\rho u ) = 0, 
\qquad \quad 
\del_t (\rho u ) + \del_x( \rho u^2 + p(\rho)) = 0,
$$
($\rho, u$ denoting the density and velocity of the fluid and $p=p(\rho)$ its pressure). The Burgers equation is also the simplest example of nonlinear hyperbolic conservation law and has played an essential role in the development of robust and accurate, shock-capturing schemes for the approximation of entropy solutions to nonlinear hyperbolic systems. 

In the present paper, we derive several {\sl relativistic generalizations of Burgers equation} which, in addition, we formulate on a {\sl curved background,} and we study the discretization of these equations via a well balanced, finite volume strategy based on time--independent solutions to these relativistic Burgers equations. More precisely, we consider the following class of {\sl nonlinear hyperbolic balance laws} 
\be
\label{AAA}
\mbox{div}_\omega (T(v)) = S(v),
\ee
posed on an $(n+1)$-dimensional curved spacetime (with boundary), denoted by $(M, \omega)$ and endowed with a volume form $\om$. The unknown function is a scalar field $v: M \to \RR$ and $\mbox{div}_\omega$ denotes the divergence operator associated with the volume form. The {\sl flux vector field} $T=T(v)$ is a field of tangent vectors prescribed on $M$ and depending smoothly upon $v$, while the right--hand side $S=S(v)$ is a scalar field prescribed on $M$ which also smoothly depends upon $v$.  

Hyperbolic conservation laws on curved spacetimes have been analyzed in recent years by LeFloch and collaborators; cf.~the review \cite{PLF}, as well as 
\cite{AmorimLeFlochOkutmustur, BL,LeFlochOkutmustur}. In order to formulate the initial value problem for (\ref{AAA}), we assume that the spacetime is foliated by hypersurfaces, that is,
\be
\label{foli}
M = \bigcup_{t \geq 0} H_t,
\ee
in which each slice $H_t$ is an $n$-dimensional manifold, canonically endowed with a normal $1$-form field $N_t$
and with the same topology as the one of the initial slice $H_0$. The {\sl global hyperbolicity} of the spacetime and the equation
(\ref{AAA}),  by definition, requires that each slice $H_t$ be {\sl spacelike,} in the
 sense that the function
\be
\label{hyperb}
v \mapsto T^0(v) := \la N_t, T(v) \ra \quad \mbox{ is strictly increasing.}
\ee
The class of nonlinear hyperbolic equations (\ref{AAA})--(\ref{hyperb}) provides us with a scalar model on which one can develop and analyze numerical methods of approximations which are then applicable to the Euler equations of relativistic compressible fluid flows on a curved spacetime. We refer the reader to the follow--up work by LeFloch and Makhlof~\cite{LeFlochMakhlof} which addresses this generalization. For the well posedness theory of the initial value problem associated with (\ref{AAA}), we refer the reader to \cite{AmorimLeFlochOkutmustur} and the references cited therein.

In this paper, we introduce a geometric formulation of the finite volume method, which is classically formulated without taking the underlying geometry into account, and we formulate a well balanced version of this method, which is inspired from the work of Russo and collaborators  \cite{PuppoRusso, Russo1,Russo2,Russo3} and is built upon the Nessyahu--Tadmor scheme \cite{NessyahuTadmor}. For a review on well balanced techniques, we refer the reader to Bouchut \cite{Bouchut}, 
Greenberg and Leroux \cite{GL}, and LeVeque \cite{Le}. 

An outline of the present paper is as follows. In Section~2, we derive our relativistic version of Burgers equation when all geometric effects are neglected. In Section~3, we explain how to take the geometry into account and we also connect our model to the Euler equations. Next, in Section~4, we present the finite volume strategy, and finally, in Section~5, several numerical experiments. Section~6 contains concluding remarks on this paper.


\section{A Lorentz--invariant conservation law} 
\label{sec-2}

\subsection{Derivation of the new model}

In this section, we search for flux vector fields $T=T(v)$ for which solutions to the equation (\ref{AAA}) satisfy the same Lorentz invariant property as the solutions to the Euler equations of relativistic fluids. For the sake of simplicity in the presentation (and without loss of generality for the purpose of this section), we assume that $n=1$, $S(v) \equiv 0$, and that the manifold $M=[0, +\infty)\times \RR$ is covered by a (single) coordinate chart $(x^0, x^1)$ in which the volume form reads  $\omega = dx^0dx^1$. Hence, (\ref{AAA}) takes the simplest form of a {\sl hyperbolic conservation law,} i.e. 
$$
\del_0 T^0(v) + \del_1 T^1(v) = 0,
$$
with $\del_0 = \del/\del x^0$ and $\del_1 = \del/\del x^1$, while $x^0 \in [0, +\infty)$ and
$x^1 \in \RR$.
In addition, in the present section,
we also assume that the functions $T^0= T^0(v)$ and $T^1= T^1(v)$ are independent of $(x^0, x^1)$.
A formulation taking geometric effects into account will be introduced and investigated in Section~\ref{sec-3}, below.

The {\sl Lorentz transformation} $(x^0,x^1) \mapsto (\xb^0,\xb^1)$ reads 
\be
\label{Lorentz}
\xb^0  := \gamma_\eps(V) \, \big( x^0 - \eps^2 V  x^1 \big),
\qquad \quad 
\xb^1 := \gamma_\eps(V) \, \big(- V \, x^0 + x^1 \big),
\ee
where $\eps \in (-1, 1)$ denotes the inverse of the (normalized) speed of light,
and $\gamma_\eps(V) = \big(1 - \eps^2 \, V^2\big)^{-1/2}$ the so-called {\sl Lorentz factor} associated with a given speed $V \in (-1/\eps, 1/\eps)$.
These equations can be inverted to give
$$
x^0  = \gamma_\eps(V) \, \big( \xb^0 +  \eps^2 V  \xb^1 \big),
\qquad \quad 
x^1  = \gamma_\eps(V) \, \big( V \, \xb^0 + \xb^1 \big).
$$
Instead of $V$, it is often convenient to use the {\sl normalized speed} $U \in \RR$ defined by
$$
e^{\eps U} := \gamma_\eps(V) \, (1 + \eps \, V) = \Bigg( {1 + \eps V \over 1 - \eps V} \Bigg)^{1/2},
$$
so that the Lorentz transformation takes the compact form
${\xb^0 \pm \eps \, \xb^1} = {e^{\mp \eps U} \, (x^0 \pm \eps \, x^1)}$
or, equivalently,
\be
\label{Lore}
\xb^0 = \cosh(\eps U) \, x^0 - \eps \, \sinh(\eps U) \, x^1,
\qquad \quad 
\xb^1  = - {\frac{1}{\eps}}\sinh(\eps U) \, x^0 +  \cosh(\eps U) \, x^1,
\ee
with $\cosh(\eps U) =  \gamma_\eps(V)$ and $\sinh(\eps U) = \eps V \, \gamma_\eps(V)$. 
Recall also that the set of Lorentz transformations (together with the spatial rotations if $n \geq 2$)
forms the so-called Lorentz isometry group, characterized by the condition that the length element of the Minkowski metric is
preserved, that is, 
$
- \eps^{-2} \, (\xb^0)^2 + (\xb^1)^2
= - \eps^{-2} \, (x^0)^2 +  (x^1)^2.
$
The relativistic Euler equations of compressible fluids are invariant under the Lorentz transformation (\ref{Lorentz}).
More precisely, given a speed $V$, the fluid velocity component $v$ in the coordinate system $(x^0, x^1)$ is related to the component $\vb$
in the coordinates $(\xb^0,\xb^1)$ by the relation
\be
\label{addition}
\vb = {v - V \over 1 - \eps^2 V\, v}.
\ee
Clearly, in the nonrelativistic limit $\eps \to 0$, one recovers the {\sl Galilean transformation}
\be
\label{Gali}
\xb^0 = x^0,
\quad
\xb^1 = - V \, x^0 + x^1, \qquad \vb = v - V  \quad \mbox{ when } \eps = 0,
\ee
in which the speed $V$ belongs to $\R$.

\begin{theorem}[Relativistic version of Burgers equation]
The conservation law
\be
\label{hf0}
\del_0 T^0(v) + \del_1 T^1(v) = 0
\ee
is invariant under the Lorentz transformations (\ref{Lore})--(\ref{addition}) if and only if,
after suitable normalization (discussed below), one has
\be
\label{hf}
T^0(v) = \frac{v}{\sqrt{1-\eps^2 v^2}}, \qquad
T^1(v) = \frac{1}{\eps^2} \Bigg( -1 + \frac{1 }{ \sqrt{1- \eps^2 v^2}}  \Bigg),
\ee
where the unknown scalar field $v$ takes its value in the interval $(-1/\eps, 1/\eps)$.
\end{theorem}

We refer to (\ref{hf0})--(\ref{hf}) as the {\bf relativistic Burgers equations} (on Minkowski spacetime).

{\sl Proof.}   We first  change  coordinates and construct certain auxilliary functions, denoted below by $H_\pm (v)$. Precisely, in the coordinates $(x^0,x^1)$, we can differentiate the Lorentz expressions (\ref{Lorentz}) and obtain 
\begin{eqnarray*}
\frac{\del \xb^0}{\del x^0} &= \frac{1}{\sqrt{1-\eps^2 \, V^2}} =\gamma_\eps(V),
\qquad \qquad \,
& \frac{\del \xb^0}{\del x^1}  =  - \eps^2 \, V\gamma_\eps(V),
\\
\frac{\del \xb^1}{\del x^0}  &\hskip-1.9cm = - V \,\gamma_\eps(V),
\qquad
& \frac{\del \xb^1}{\del x^1} = \gamma_\eps(V) 
\end{eqnarray*} 
and, by the chain rule, 
\begin{eqnarray*}
 \del_{0} T^0
& =&\frac{\del T^0 }{\del \xb^0}\frac{\del \xb^0 }{\del x^0}
+ \frac{\del T^0 }{\del \xb^1}\frac{\del \xb^1 }{\del x^0}
=  \frac{\del T^0}{\del \xb^0} \gamma_\eps(V)
 - \frac{\del T^0}{\del \xb^1}\gamma_\eps(V) \,V,
\\
 \del_{1} T^1
& =& \frac{\del  T^1}{\del \xb^0}\frac{\del \xb^0 }{\del x^1}
+ \frac{\del T^1 }{\del \xb^1}\frac{\del \xb^1 }{\del x^1}
=  - \eps^2 \, \frac{\del  T^1}{\del \xb^0} \gamma_\eps(V) \,V+  \frac{\del T^1 }{\del \xb^1} \gamma_\eps(V).
\end{eqnarray*} 
By substituting these formulas in (\ref{hf0}), the conservation law reads 
\begin{eqnarray*}
0=\del_{0} T^0 + \del_{1} T^1
=\del_{0} (\gamma_\eps(V) \, T^0  -  \eps^2 \, \gamma_\eps(V) \, V \, T^1) + \del_{1}(- V\,\gamma_\eps(V) \,T^0
+ \gamma_\eps(V) \,T^1),
\end{eqnarray*}
which has the {\sl same structure} as (\ref{hf0}) in the coordinates $(\xb^0,\xb^1)$  {\sl provided}
\begin{eqnarray}
\label{EX2} 
\gamma_\eps(V) \,\big(T^0(v) - \eps^2V \, T^1(v) \big) &= T^0\Bigg(\frac{v-V}{1-\eps^2v\,V}\Bigg) + C_1,
\\
\gamma_\eps(V) \,\big(- V T^0(v) + T^1(v)\big)&
=T^1\Bigg( \frac{v-V}{1-\eps^2v\,V}\Bigg)+C_2,
\end{eqnarray}
for some constants $C_1, C_2$. This leads us to the desired invariance property and it remains to determine the general expressions of $T^0, T^1$.

We can always normalize the flux terms so  that $T^0(0)=T^1(0)=0,$ and then observe that the constants $C_1=-T^0(-V)$ and $ C_2=-T^1(-V)$ are determined as functions of $V$. Next, multiplying the second equation in (\ref{EX2}) by $\pm \eps$ and summing up with the first one, we find  
$$
\gamma_\eps(V) \, (1\mp\eps V)(T^0\pm \eps \, T^1)(v)= \big( T^0\pm \eps \, T^1 \big) \, \Bigg(\frac{v-V}{1-\eps^2\, v V}\Bigg) - (T^0\pm \eps \, T^1)(-V).
$$
By setting $(T^0 \pm \eps \, T^1)= : H_{\pm}$, these two equations read
\be
\label{E3}
H_{\pm}(v)=
\sqrt{\frac{1\pm\eps V}{1 \mp \eps V}}\, \Bigg(H_{\pm}\Big(\frac{v-V}{1-\eps^2\, v V}\Big)
- H_{\pm}(-V)\Bigg).
\ee

Next, we define the new variables 
$$
u := \frac{1}{2\eps} \ln\Big({\frac{1+ \eps v}{1- \eps v}}\Big), \qquad \phi_{\eps}(u): =v, \quad \phi_{\eps}(U):= V,
$$
so that
$$
v= \frac{1}{\eps}\frac{e^{2\eps u} - 1}{e^{2\eps u} + 1}= \phi_{\eps}(u),
\qquad
\vb = \frac{1}{\eps} \frac{e^{2\eps (u- U) }- 1}{e^{2\eps (u-U) } +1}
=\phi_{\eps}(u-U),
\qquad
\sqrt{\frac{1\pm \eps V}{1 \mp \eps V}  }= e^{\pm \, \eps U}.
$$
By inserting these expressions in (\ref{E3}), it follows that
\be
\label{E3'}
H_\pm \big( \phi_{\eps}(u)\big) -H_\pm \big( \phi_{\eps}(0)\big)  =  e^{\pm \, \eps U}
\Big( H_\pm \big( \phi_{\eps}(u- U)\big) - H_\pm \big( \phi_{\eps}(-U)\big) \Big).
\ee

Abusing notation, we write $H,\phi, e^{\eps U}$ instead of $H_\pm, \phi_\eps, e^{\pm \eps U}$,  respectively. We differentiate (\ref{E3'})  with respect to $u$ and obtain 
$H' \big( \phi(u)\big) \phi'(u)     = e^{ \eps U}  H' \big( \phi(u- U)\big) \phi'(u- U)$ 
and, next, with respect to $U$
$$
H''\big( \phi(u-U) \big) \big(\phi'(u-U)\big)^2 = \eps  H'\big(\phi (u-U)\big) \phi'(u-U)  - \phi''(u-U) H'\big( \phi(u-U) \big).
$$
By letting $U=0$, it follows that
$H''\big( \phi(u) \big) \big(\phi'(u)\big)^2 = \eps  H'\big(\phi (u)\big) \phi'(u)  - \phi''(u) H'\big( \phi(u) \big)$.
The substitution  $ K(u): = H(\phi(u))$  allows us to rewrite the latter equation as $K''(u) = \eps \, K'(u)$, 
which has the following general solution (in which $c,  \widetilde{c}$ are constants): 
$$
K(u) = c \frac{e^{\eps u}}{\eps} + \widetilde{c}.
$$
Alternatively, this latter relation can be derived by dividing (\ref{E3'}) by $u \neq 0$ and letting $u \to 0$. 

Returning to $H_\pm$, we deduce that there exist constants so that 
$H_\pm \big(\phi(u)\big) = \pm \eps^{-1} \big( e^{\pm \eps u}- 1 \big)$. 
Recalling that $H_\pm = T^0 \pm \eps T^1$  with $H_\pm\big(\phi(u)\big)= H_\pm(v)$, we find
$
T^0= \frac{H_+ + H_-}{2 }$ and $T^1= \frac{H_+- H_-}{2\eps},
$
from which we  deduce
\begin{eqnarray*}
T^0(v)&=& T^0\big(\phi(u)\big)= \frac{e^{\eps u} - e^{-\eps u}}{2\eps} =  \frac{1}{\eps} \sinh(\eps u)
= u + O ( \eps ^2u ^3), \\
T^1(v)&=& T^1\big(\phi(u)\big)= \frac{ e^{\eps u} + e^{-\eps u} -2 }{2 \eps^2} = \frac{1}{\eps^2} \big( \cosh (\eps u) -1  \big)
= \frac{u^2}{2} + O ( \eps ^2u ^4).
\end{eqnarray*}
Observe that $T^0$ and $T^1$ are linear and quadratic, respectively, and that, in the limiting case $\eps\to 0$,
we recover the inviscid Burgers equation. Moreover, substituting
$u= \frac{1}{2\eps} \ln\big({\frac{1+\eps v}{1- \eps v}}\big)$ in $T^0\big(\phi(u)\big)$ and $T^1\big(\phi(u)\big),$
we thus arrive at (\ref{hf}).  $\Box$


\subsection{Hyperbolicity and convexity properties}

The proposed relativistic Burgers equation retains several key features of the relativistic Euler equations:
\begin{itemize}

\item Like the Euler system, it has a conservative form.

\item Like the velocity component in the Euler system, our unknown $v$ is constrained to lie in the interval
$(-1/\eps,1/\eps)$  limited by the (inverse of the) light speed parameter.

\item Like for the Euler system, by sending the light speed to infinity one recovers the classical (nonrelativistic) model.
\end{itemize}

\begin{theorem}[Properties of the relativistic Burgers equation]\
\begin{enumerate}

\item The map $w= T^0(v) = \frac{v}{\sqrt{1- \eps^2 v^2}}\in \RR$ is increasing and  one-to-one from $(-1/\eps, 1/\eps)$ onto $\RR$.

\item In terms of the new unknown $w\in \RR$, the equation (\ref{hf0}) is equivalent to
\be
\label{wf}
\del_{0}  w + \del_{1} f_\eps (w)= 0,
\qquad \quad 
f_\eps(w)=  \frac{1 }{\eps^2} \Big( -1 +  \sqrt{1+ \eps^2 w^2}  \,\Big).
\ee
The flux  $f_\eps$ is strictly convex and, therefore, the conservation law (\ref{hf0}) is genuinely nonlinear in the sense that
$
\frac{\del_v T^1 (v)}{\del_v T^0(v)}= \del_w {f_\eps}(w)
$
is strictly increasing or strictly decreasing in the variable $T^0(v)$.  

\item In the nonrelativistic limit $\eps\to 0$, one recovers the inviscid Burgers equation
\be
\label{Burg} 
\del_0 u + \del_1 ( u^2/2 )=0, \qquad u= u(x^0, x^1) \in \RR. 
\ee  
\end{enumerate}
\end{theorem}

{\sl Proof.} We set  $w:= T^0(v) = \frac{v}{\sqrt{1- \eps^2 v^2}}\in \RR$ and take its derivative  with respect to $v$
$$
\del_v T^0(v) = \frac{1}{(1- \eps^2 v^2)^{3/2}} > 0, \qquad v\in {(-1/\eps, 1/\eps)},
$$
so that $T^0(v)$ is increasing and  one-to-one from $(-1/\eps, 1/\eps)$ onto $\RR$. We get 
$v= \frac{w}{\sqrt{1 + \eps^2 w^2}}$ (after choosing the positive branch, for the sake of definiteness) and, substituting in (\ref{hf0}),  we arrive at (\ref{wf}). The flux is strictly convex, since 
$$
f_\eps'(w)  =\frac{w}{\sqrt{1+ \eps^2 w^2}}, \qquad
f_\eps''(w) = \frac{1}{({1+ \eps^2w^2})^{3/2}} > 0.
$$ 
Furthermore, we have
$$
\frac{\del_v T^1 (v)}{\del_v T^0(v)}= v = \frac{w}{\sqrt{1+ \eps^2 w^2}} = f_\eps'(w),
$$
which is positive for $w>0$ and negative for $w<0.$ Finally, we have  
$$
T^0(v)=  \frac{1}{\eps} \sinh(\eps u)= u + O ( \eps ^2u ^3), \qquad
T^1(v)=  \frac{1}{\eps^2} \big( \cosh (\eps u) -1  \big)= \frac{u^2}{2} + O ( \eps ^2u ^4).
$$
When $\eps$ approaches $0$, one recovers the standard inviscid Burgers equation (\ref{Burg}).
$\Box$


\subsection{The nonrelativistic limit}

The Galilean transformation $(x^0,x^1) \mapsto (\xb^0,\xb^1)$ is defined by (\ref{Gali}), in which 
the speed $V \in \RR$ is a given parameter. The velocity component $v$ in the coordinate system $(x^0,x^1)$ is related to the velocity $\vb$ in the coordinates $(\xb^0, \xb^1)$  by the relation $\vb=  v - V$. Recall that one recovers the Galilean transformations from the Lorentz transformations in the limit  $\eps \to 0$.

\begin{theorem} [A derivation of the (nonrelativistic) Burgers equation]
The conservation law
\be
\label{GAL2}
\del_0T^0(v) + \del_1 T^1(v)= 0
\ee
is invariant under the Galilean transformations (\ref{addition})--(\ref{Gali}) if and only if the flux $T^0$ and $T^1$ are linear and quadratic functions, respectively. After a suitable normalization, one has $T^0(v)=v$ and $T^1(v)= v^2/2$.
\end{theorem}

{\sl Proof.} Using the chain rule, we find 
$$
\del_0 T^0= \frac{\del T^0}{\del \xb^0}\frac{\del \xb^0 }{\del x^0} + \frac{\del T^0 }{\del \xb^1 } \frac{\del \xb^1}{\del x^0}
= \frac{\del T^0}{\del \xb^0}  - V \frac{ \del T^0}{ \del \xb^1},
\qquad\qquad
 \del_1 T^1 = \frac{\del T^1}{\del{ \xb^1}}.
$$
By substituting this in (\ref{GAL2}), it follows that
$$
\del_0 T^0(\vb + V) + \del_1 \Big(T^1(\vb + V) - V T^0 (\vb +V)\Big)= 0.
$$
Now, this equation has the same algebraic structure as (\ref{GAL2}), provided
\be
\label{IN1}
T^0(\vb+ V)  = T^0(\vb) + C_1,
\qquad \qquad 
T^1(\vb +V) - V\, T^0(\vb + V) = T^1(\vb) + C_2,
\ee
where $C_1 $ and $ C_2$ are constants depending upon $V$.

Then, we can derive the first equation in (\ref{IN1}) with respect to $\vb$, that is, 
$
T^0_{\vb}(\vb+V) = T^0_{\vb}(\vb),
$
from which we conclude that $T^0_{\vb}$ is periodic of period $V$. As a result,   $T^0$ is a linear function of the form $T^0(\vb) = C\, \vb + \widetilde{C}$ for some constants $C$ and $\widetilde{C}$.
Next, multiplying the first equation in (\ref{IN1}) by $V$ and adding up to the second one, we get
$
T^1(\vb +V) =  T^1(\vb) +   V\, T^0(\vb) + V C_1 + C_2.
$
Since $T^0$ is a linear function and deriving  twice the given equation,
 it follows that
\begin{eqnarray*}
T^1(\vb +V) - T^1(\vb) & = & C V\, \vb +V \widetilde{C}+  V C_1 + C_2, 
\\
T^1_{\vb}(\vb + V) -  T^1_{\vb}(\vb) - C V & =& 0, 
\\
T^1_{\vb\,\vb}(\vb + V) - T^1_{\vb\,\vb} (\vb)& =& 0.
\end{eqnarray*}
Therefore,  $T^1_{\vb\, \vb}$ is periodic with period $V$, which implies that the first--order derivative $T^1_{\vb} $ is linear and that  $T^1$ is quadratic. After normalization,  one obtains exactly Burgers equation,  corresponding $T^0(v)= v$ and $T^1(v)= v^2/2$. 
$\Box$ 


\subsection{Derivation from the relativistic Euler equations}

\subsubsection*{Constant density}

We present another derivation of (\ref{wf}), now from the relativistic Euler equations which read 
\begin{eqnarray}
\label{REE}
&\del_0 \bigg( \frac{p(\rho)+\rho c^2}{c^2}\frac{v^2}{c^2 - v^2} + \rho\bigg) + \del_1 \bigg( (p(\rho) + \rho c^2) \frac{v}{c^2 - v^2}\bigg) = 0,
\\
&\del_0\bigg( (p+ \rho c^2) \frac{v}{c^2 - v^2}\bigg) + \del_1\bigg( (p(\rho) +\rho c^2) \frac{v^2}{c^2 - v^2} +p(\rho) \bigg)=0,
\nonumber 
\end{eqnarray} 
where $\rho$ and $u$ denote the mass--energy density and the (suitably normalized) velocity of the fluid, respectively, while the pressure $p=p(\rho)$ is given by an equation of state depending on the fluid and the constant $c>0$ represents the light speed.  

By formally setting the density $\rho$ (and thus the pressure $p(\rho)$) to be a constant in the second Euler equation above, we obtain 
$$
\del_0\bigg( \frac{v}{c^2 - v^2}\bigg) + \del_1\bigg( \frac{v^2}{c^2 - v^2} \bigg)=0.
$$
Introducing the change of unknown $z= \frac{v}{1-\eps^2v^2}$ or $v = \frac{-1 \pm \sqrt{1+ 4 \eps^2 \,z^2 }}{2\eps^2\, z}$, with 
 $c=1/\eps$, we find 
$$
\del_0 z  + {1 \over 2\, \eps^2}  \del_1 \bigg(-1 \pm \sqrt{1+ 4\eps^2 z^2}\bigg) = 0,
$$
which is precisely (\ref{wf}) if we choose the positive branch and suitably normalize the equation. 


\subsection{Derivation from a general covariant equation on Minkowski spacetime}
\label{sec:25}

We now adopt a more geometric standpoint, but yet we begin by the case that the metric is flat. Following the notation in the introduction, we consider a general conservation law posed on a spacetime $(M,g)$ endowed with a Lorentzian metric $g$, i.e. 
\be
\nabla_\alpha T^\alpha(v) =0 \quad \mbox { in }   M,
\ee
where $v: M \to  \RR$ is the unknown scalar field, and $T^\alpha =T^\alpha(v)$ is a prescribed vector field depending on $v$ as a parameter. 
By analogy with the Euler equations, we impose that, for every constant $v$, $T^\alpha(v)$ is a unit vector field, that is, 
$g_{\alpha \beta} T^\alpha (v) T^\beta (v) = 1$. 

Specifically, we now assume that $(M,g)$ is the $(1+1)$--Minkowski spacetime, described in coordinates $(x^0,x^1)= (t,x)$, 
by $g=-\eps^{-2} dt^2 + dx^2$. We obtain $-\eps^{-2} T^0(v)^2 + T^1(v)^2 = 1$, 
thus $T^1(v) = \pm \sqrt{1+ \eps^{-2} T^0(v)^2}$ and, after choosing the positive branch,  
$$
\del_0 T^0 + \del_1 \big(C + \sqrt{1+ \eps^{-2} (T^0)^2} \big) = 0. 
$$
Here,  the constant $C$ was added for convenience. Finally, introducing the change of variable $w= \eps^{-2} T^0$ and normalizing the constant to be $C=-1$, we again recover the proposed relativistic Burgers equation (\ref{wf}). 


\section{Burgers equations with geometric effects}
\label{sec-3}

\subsection{Relativistic Euler equations on Schwarzschild spacetime}
 
\subsubsection*{Vanishing pressure on flat spacetime} 

Our standpoint is now different and we start directly from the Euler equations and formally assume that, in the relativistic Euler equations (\ref{REE}), the pressure vanishes identically, that is, 
$$
\del_0\Bigg( \frac{\rho}{c^2 - v^2} \Bigg) + \del_1 \Bigg( \frac{\rho \,v }{ c^2 - v^2}\Bigg) = 0, 
\qquad \quad 
\del_0\Bigg( \frac{\rho \,v}{c^2 - v^2}\Bigg)  + \del_1 \Bigg( \frac{\rho \,v^2 }{ c^2 - v^2}\Bigg) = 0.
$$
These two equations can be rewritten as 
\begin{eqnarray*}
&(c^2-v^2)(\del_0\rho + v \del_1 \rho) + \rho( 2v \del_0v + (v^2 + c^2)\del_1v)=0, 
\\
&v(c^2-v^2)(\del_0\rho + v \del_1 \rho) + \rho( (v^2 + c^2)\del_0v + 2vc^2 \del_1v)=0,
\end{eqnarray*} 
and we recover the {\bf classical Burgers equation} (in flat spacetime) 
\be
\label{standard}
\del_0 v + \del_1(v^2/2) =0,
\ee
which, therefore, appears to be also a relevant model, even in the relativistic setting.

 
\subsubsection*{Schwarzschild geometry}

We now want to take geometric effects into account in (\ref{standard}). The scalar equation of interest is still sought in the form 
\be
\nabla_\alpha T^\alpha(v) = 0, \qquad v: M \to \RR
\label{1.1}
\ee
for some prescribed vector field $T^\alpha(v)$. To determine the model of interest, we start from the Euler equations posed on a curved spacetime and, for definiteness, assume that the background manifold $(M,g)$ is Schwarzschild spacetime. The latter geometry 
describes a spherically symmetric black hole solution of the vacuum Einstein equations given, in suitably chosen coordinates
$(x^0, x^1, x^2, x^3) = (c t,r, \theta, \varphi)$ adapted to the symmetry, by 
\be
\label{metric}
g = -\biggl( 1 - {2m \over r} \biggr) c^2 \, dt^2
  + \biggl( 1 - {2m \over r} \biggr)^{-1} \, dr^2
  + r^2 \bigl( d \theta^2 + \sin ^2 \theta \, d\varphi^2 \bigr)
\ee
for $t>0$ and $r > 2 m$ where the coefficient $m>0$ is referred to as the mass. 
This spacetime is {\sl spherically symmetric,} that is invariant under the group of
rotations operating on the space-like $2$-spheres given by keeping $t$ and $r$ constant.
It is {\sl static,} since the vector field $\del_t$ is a time-like Killing vector, and is asymptotic to (flat) Minkowski in the limit $r \to \infty$.
In (\ref{metric}), there is an apparent (not physical) singularity at $r= 2m$ which, however, is removed in the so--called Eddington-Finkelstein coordinates. 

The Euler equations read 
\be
\label{Euler0} 
\nabla_\alpha T^\alpha_\beta = 0, \qquad \quad T_\alpha^\beta = (\rho c^2 + p) \, u_\alpha \, u^\beta + p \, g_\alpha^\beta,
\ee
where $T^\alpha_\beta$ is the energy-momentum tensor of perfect fluids, 
$\rho$ is the mass--energy density of the fluid, and $u=(u^\alpha)$ denotes its unit velocity vector, so $g_{\alpha\beta} \, u^\alpha u^\beta = -1$. 
This model is supplemented by an equation of state for the pressure $p = p(\rho)$. 
We assume that solutions to the Euler equations depend only on the time
variable $t$ and radial variable $r$, and that the non-radial components of the velocity vanish, that is, 
$(u^\alpha) = (u^0(t,r), u^1(t,r), 0,0)$. Since $u$ is unit vector, we obtain
$$
-1 = - \Biggl(1 - {2m \over r}\Biggr) \, (u^0)^2 + {(u^1)^2 \over \Bigl(1 - {2m \over r}\Bigr)}.
$$
It is convenient to introduce the velocity component 
$v :=  c \bigl(1 - {2m \over r}\bigr)^{-1} \, {u^1 \over u^0}$, 
so that
$$
(u^0)^2 = {c^2 \over (c^2-v^2) \bigl(1 - {2m \over r}\bigr)}, \qquad
\quad 
(u^1)^2 =\Bigl(1 - {2m \over r}\Bigr) { v^2 \over (c^2-v^2)}.
$$
We obtain 
$$
T^{00} = {\tT^{00} \over \bigl( 1 - {2m \over r} \bigr)}, \quad
T^{01} = \tT^{01}, \quad
T^{11} =  (1 - {2m \over r}) \, \tT^{11}, \qquad T^{22} := {p(\rho) \over r^2}, 
$$
with
$$
\tT^{00} := {c^2 \rho + p(\rho) v^2/c^2 \over c^2 - v^2} c^2, 
\quad
 \tT^{01} := {c^2 \rho + p(\rho) \over c^2 - v^2} cv, 
\quad
 \tT^{11} := {v^2 \rho + p(\rho) \over c^2 - v^2} c^2.
$$ 


\subsubsection*{Vanishing pressure on Schwarzschild spacetime}

We suppose first that the pressure $p$ vanishes identically, and the Euler system in $1+1$ dimensions on a Schwarzschild background takes the simplified form: 
\begin{eqnarray*}
& \del_t \Bigg( \frac{r^2}{c^2}\, \tT^{00} \Bigg) + \del_r \Bigg(\frac{r(r-2m)}{c} \, \tT^{01}\Bigg) = 0,
\\
& \del_t \Bigg(\frac{r(r-2m)}{c} \,\tT^{01}\Bigg)
 + \del_r \Big((r - 2m)^2 \, \tT^{11} \Big) - 3m \, \frac{(r-2m)}{r} \, \tT^{11}
 + m \, \frac{(r-2m)}{r}  \tT^{00}=0 
\end{eqnarray*} 
or, equivalently,
\begin{eqnarray}
\label{EE1}
& \del_t \Bigg( \frac{r^2 \rho}{c^2-v^2} \Bigg)
+ r^2 \Big(1-\frac{2m}{r} \Big)\del_r \Bigg( \frac{\rho \, v}{c^2-v^2}\Bigg) + \frac{2r-2m}{c^2-v^2} \,\rho v= 0,
\\
& \del_t \Bigg( \frac{\rho\, v}{c^2-v^2} \Bigg)
 + \Big( 1-\frac{2m}{r} \Big) \del_r \Bigg( \frac{v^2 \rho}{c^2 -v^2}\Bigg) + \rho\, \frac{m(c^2 -3v^2)+ 2v^2r}{r^2 (c^2-v^2)}
 = 0.
\label{EE1-2}
\end{eqnarray} 
By combining these two equations together, we arrive at 
$$
\del_t v + \Big( 1 - \frac{2m}{r} \Big) \, v\del_r v = \frac{m}{r^2}(v^2-c^2) 
$$ 
or, equivalently,
\be
\label{RBE}
\del_t(r^2\,v) + \del_r \Big( r ( r - 2m) \, \frac{v^2}{2} \Big) = rv^2 - m c^2. 
\ee
This provides us with a variant of Burgers equation which incorporates two features of Schwarz\-schild spacetime, i.e. the effect of spherical symmetry and the mass of the black hole. We refer to it as {\bf Burgers equation on Schwarzschild spacetime}.


\subsection{Static solutions of Burgers equation on Schwarzschild spacetime}

\subsubsection*{Static solutions} 

Static solutions will play a central role in order to design our numerical method. 
We search for $t$--independent solutions to the equation (\ref{RBE})  
\be
\label{SS3}
 {1 \over 2}  \del_r \Big( r ( r -2m) \, v^2 \Big) = rv^2 - m c^2.
\ee
Using the change of variable $X:=c^2-v^2$ and $Y:= (r-m)/{m}$, we find 
$\frac{dX}{X}= -2\, \frac{dY}{1- Y^2}$ and, after integration, 
$
X= K^2\, \frac{Y-1}{Y+1},
$
where $K \in (0,c)$ is an arbitrary constant. We can thus describe all steady solutions $v=v_s(r)$ by
\be
\label{SS2-00}
v_s(r) =\pm \sqrt{c^2 - K^2 \Big(1-\frac{2m}{r} \Big)}, \qquad \quad K \in (0,c)
\ee
or, equivalently,
\be
\label{SS2}
{c^2 - v_s^2 \over 1- \frac{2m}{r}} = K^2, \qquad \quad K \in (0,c).
\ee

\subsubsection*{Conservative form}

In fact, the above calculations suggest to check whether the equation (\ref{RBE}) admits the conservative form 
$$
\del_0(\alpha\,v) + {1 \over 2} \del_y \Big( \alpha\, \big( v^2-c^2 \big) \Big)=0,
$$
and we find that $y=y(r)$ and $\alpha= \alpha(y)$ must satisfy
$$
\frac{\alpha'(y)}{\alpha(y)}= \frac{-2m}{r^2},
\qquad \quad dy=\frac{dr}{1-\frac{2m}{r}}.
$$
Solving these equations  yields $\alpha=\frac{r}{r-2m}$ and $y= r + 2m\, \ln\,(r-2m)$, and we obtain the following {\bf conservative form of  Burgers equation on Schwarzschild spacetime}
\be
\label{CF3}
\del_t \Big(\frac{r}{r-2m}\,v\Big) +  {1 \over 2}  \del_r \Big( \frac {r}{r-2m} \, \big( v^2-c^2 \big) \Big)=0.
\ee
Interestingly the variable $y$ above is essentially the Eddington--Finkelstein coordinate.


\subsection{Relativistic Burgers equation with geometric effects}

We now return to the approach adopted in Section~2.  
We assume that the manifold is described by a single chart and, after identification, we set
$M= \RR_+ \times \RR$. In coordinates $(x^0, x^1)$ with $\del_\alpha := \del/\del x^\alpha$
(with $\alpha= 0,1$), the hyperbolic balance law under consideration reads 
\be
\label{CL.1}
\del_0 \big( \omb \, T^0(v, x^0, x^1) \big) + \del_1 \big( \omb \, T^1(v,x^0, x^1) \big) = \omb \, S(v, x^0, x^1),
\ee
where $v: M \to \RR$ is the unknown function and $T^\alpha = T^\alpha(v, x^0, x^1)$ and $S=S(v, x^0, x^1)$ are prescribed (flux and source) fields on $M$, while $\omb =\omb(x^0, x^1)$ is a positive weight-function.
This equation is {\sl hyperbolic in the direction} $\del/\del x^1$ provided
$\del_v T^0 (v) >0$, which is always assumed.

For instance, in the case that the flux and source are independent of $x^0, x^1$, the above balance law can be rewritten in the form
\begin{eqnarray*}
& \del_0 v + \del_1 f(v) = \St(v),
\qquad
\quad \del_v f(v)
:= { \del_v T^1(v) \over \del_v T^0(v)}, 
\\
&\St(v)
:=  {1 \over \del_v T^0(v)} \, \Big( S(v) - \del_0 \Omb \, T^0(v) - \del_1 \Omb \, T^1(v) \Big), \qquad \Omb := \ln \omb.
\end{eqnarray*} 
It may occur that $\St$ vanishes identically,  however, one should keep in mind that (\ref{CL.1}) is the {\sl geometric form} of this equation: it should  preferred for the numerical discretization. It may also occur that $S$ vanishes identically, in which case the original equation (\ref{CL.1}) is a conservation law, but the expanded form is nonconservative. 


More precisely, within the above setting, we now take into account geometric effects in the relativistic Burgers equation (\ref{wf}). 
By analogy with the first part of this section and in view of the equation (\ref{RBE}), we set $x^0=c t$ and $x^1 =r$, with $c=1/\eps$
and propose the following model:  
\begin{eqnarray}
\label{ss} 
&\del_t (r^2 \, w) + \del_r \Big(r(r-2m) \, f_\eps(w) \Big) = 0, 
\qquad
\quad
f_\eps(w)= \frac{1}{\eps^2} \Big( -1 + \sqrt{ 1+ \eps^2 w^2} \, \Big),
\end{eqnarray} 
which we refer to as the {\bf relativistic Burgers equation on Schwarzschild spacetime.} The choice of flux and source-term in (\ref{CL.1}) made here 
allows us to establish a clear connection with the model (\ref{RBE}) derived in the previous subsection from the physicaly sound assumption of Lorentz invariance. 

Next, let us seek for steady solutions $w_s = w_s(r)$, satisfying 
$$
\del_r \Big( r(r -2m) \, f_\eps(w) \Big) = 0
$$
and recalling that $f_\eps$ is strictly convex, we can define its inverse $f_\eps^{-1}$ by selecting the ``positive branch'', and we find 
$w_s(r) = f_\eps^{-1}\Big( \frac{K}{r(r-2m)} \Big)$
where $K$ is constant, so that all steady solution to (\ref{ss}) are given by 
\be
\label{static5} 
w_s (r) = \pm \frac{K}{r(r-2m)} \sqrt{  \eps^2 +2{r(r-2m)\over K}}. 
\ee
  

\subsection{Summary} 

We have thus derived two Burgers--type models that appear to be of  particular interest in the context of relativistic fluid dynamics on a curved spacetime. On one hand, we have proposed the Burgers equation on Schwarzschild spacetime (\ref{RBE}), which depends on two parameters $m,c$ and reduces to the standard Burgers equation in the limit of zero mass $m$ and infinite light speed $c$. 
On the other hand, motivated by our derivation of a Lorentz invariant model described in the previous section, we have proposed the 
relativistic Burgers equation with geometric effects (\ref{ss}), which depends on the light speed $1/\eps$ and a geometric weight function $\omb$. 
Both of these models will be studied in the numerical section below.

\subsection{Alternative approaches}

\subsubsection*{Timelike flux vector}

Observe that, alternatively, we can normalize $T^\alpha(v)$ to be unit timelike, 
$g_{\alpha \beta} T^\alpha (v) T^\beta (v) = -1$, so that  $T^1(v) = \pm \sqrt{-1+ \eps^{-2} T^0(v)^2}$.  So, after normalization, we find
\be
\label{new3}
\del_0 T^0 + \del_1 \sqrt{-1+ \eps^{-2} (T^0)^2} = 0. 
\ee


\subsubsection*{Derivation from a general covariant equation on Schwarzschild spacetime}

We return to the analysis in Section~\ref{sec:25} but now work with the Schwarzschild  metric, that is, 
$$
g= -\eps^{-2} (1-2m/r) \, dt^2 + (1-2m/r)^{-1} \, dr^2
$$
with $(x^0, x^1) = (t, r)$. 
We require that the flux vector field $(T^0(v), T^1(v))$ is of unit norm, i.e.  
$$
 -\epsilon^{-2}  (1-2m/r) \, (T^0(v))^2  + (1-2m/r)^{-1} \,  (T^1(v))^2 =1,
$$
which yields 
$$
(T^1)^2 = \Bigg( 1 - {2m \over r} \Bigg) \, \Bigg( 1+ \epsilon^{-2} \Bigg(1 - {2m \over r} \Bigg) (T^0)^2 \Bigg).
$$
Finally, in term of $w: = \eps^{-2} T^0$, the associated equation can be written as 
\be
\label{2Mod}
\del_0 w + {1 \over \eps^2} \del_1\Bigg(  -1 +  \Bigg( 1 - {2m \over r} \Bigg)^{1/2} \, \Bigg( 1+ \eps^2 \Bigg(1 - {2m \over r} \Bigg) w^2 \Bigg)^{1/2}\Bigg)=0,
\ee
which, in the limit $m \to 0$, reduces to the model (\ref{wf}) derived in flat space. Then, we can introduce 
$z^2 := (1-2m/r) w^2$ and arrive at 
\be
\label{2Mod-3}
\del_0\Bigg( \Bigg( 1 - {2m \over r} \Bigg)^{-1/2} \, z \Bigg)  + {1 \over \eps^2} \del_1\Bigg(  -1 +  \Bigg( 1 - {2m \over r} \Bigg)^{1/2} \, \sqrt{ 1+ \eps^2 z^2}\Bigg)=0. 
\ee


\section{Well--balanced finite volume approximations with geometric effects}

\subsection{Geometric formulation of finite volume schemes}

We follow Amorim, LeFloch, and Okutmustur \cite{AmorimLeFlochOkutmustur} and formulate the finite volume scheme for a hyperbolic balance law posed on a curved spacetime $(M, \omega)$, i.e. 
\be
\label{HBL}
\mbox{div}_\omega (T(v))
=
{1 \over \om} \, \Big( \del_0(\om\,T^0(v, x)) + \del_1(\om \, T^1(v, x)) \Big)
= S(v, x),
\ee
where $v: M\to \RR$ is the unknown function and $T(v) = (T^\alpha(v))$ 
and $S(v): M \to \RR$ are prescribed flux and source fields, respectively, while $\om=\om(x) $ is a positive weight function which is identified with the volume form $\omega \, dx$ prescribed in some (fixed) global chart of coordinates. More precisely, we work on a $(1+1)$-dimensional spacetime and assume the global hyperbolicity condition mentioned in the introduction. 

The finite volume methodology applies to general triangulations $\Tcal^h = \bigcup_{K\in \Tcal^h} K$
of the manifold $M$ which are made of spacetime elements $K \subset M$ satisfying the following structural conditions:
\begin{itemize}
\item The boundary $\dK$ of an element $K$ is  piecewise smooth
$\dK = \bigcup_{e\subset \dK} e$ and  contains exactly two
spacelike faces (in the sense (\ref{hyperb})),
denoted by  $\ekp$ and $\ekm$, and ``timelike'' elements
$$
e^0 \in \del^0 K:= \del K \setminus \big\{\ekp, \ekm\big\}.
$$
\item The intersection $K \cap K'$ of two distinct elements $K, K'$
is  a common face of $K, K'$.
\end{itemize}
It will be convenient to denote by
 $|K|$, $|\ekp|$, $|\ekm|$, $|e^0|$ represent the (one-dimensional or two-dimensional)
measures of $K$, $\ekp, \ekm,e^0$, respectively.

We define the finite volume approximations by formally averaging the balanced  law (\ref{HBL}) over each element $K \in  \Tcal^h$ of the triangulation. By integrating in space and time, we can write
$$
\int_K \mbox{div}_\omega (T(v))\,  \om =  \int_K S(v) \,  \omega, 
$$ 
where $\omega=\omega(\cdot, \cdot)$ is here regarded as a two-form in $M$. 
With the Stokes formula, we find 
$$
\int_{\ekp}T^0(v)\, \om(n_{\ekp},\cdot)  =  \int_{\ekm}  T^0(v)\,  \om(n_{\ekm},\cdot) 
 - \sum_{e^0  \in \del^0 K}\int_{e^0} T^1(v)\,\om(n_{e^0},\cdot) +
 \int_{ K} S(v)\,\om.
$$
Given averages $v^-_K$ and $v^-_{K_{e^0}} $ computed along $\ekm$ and $e^0\in \del^0 K$,
we need to compute the average $v^+_K$ along $\ekp$. To this end, we introduce 
$$
\int_{\ekm}  T^0(v)\, \om(n_{\ekp},\cdot) \simeq  |\ekm| \,\overline{\om}_{\ekm}\, \overline{T}_{\ekm}(v_K^-),
\qquad 
\quad \int_{e^0} T^1(v)\,\om(n_{e^0},\cdot) \simeq |e^0|\, \overline{\om}_{e^0}\, Q_{K,e^0}( v_K^-, v_{K_{e^0}}^-),
$$
and
$$
\int_{ K} S(v)\,\om  \simeq |K| \, \om_K\,  \overline{S}_{ K},
$$ 
in which, for each $K$ and $e^0 \in \del ^0 K$ we have selected a (Lipschitz continuous) numerical flux 
$Q_{K,e^0} : \RR^2 \to \RR$ satisfying natural properties of consistency, conservation, and monotonicity. 
In turn, the finite volume method of interest takes the form 
\be
\label{SHM} 
\qquad 
|\ekp| \,\overline{\om}_{\ekp} \,\overline{T}_{\ekp}(v_K^+)=
|\ekm| \,\overline{\om}_{\ekm}\, \overline{T}_{\ekm}(v_K^-) - \sum_{e^0\in \del ^0 K} |e^0| \, \overline{\om}_{e^0}  \,Q_{K,e^0} ( v_K^-,v_{K_{e^0}}^-)
+    |K|\, \om_K \, \overline{S}_{ K}.
\ee
For the sake of stability of this scheme, a standard CFL condition is assumed. See  \cite{AmorimLeFlochOkutmustur} for further details. 


\subsection{Finite volume schemes in coordinates} 

From now on, motivated by the models discussed in the previous two sections, we assume that the spacetime is described in coordinates $(t,r)$ such that the weight $\omega=\omega(r)$ as well as the flux and source terms depend only upon the spatial variable.  
We divide into equally spaced cells $I_j= [ \rjmundemi,\rjpundemi]$ of size $\Delta r$, centered at $\rjjj$, 
that is,
$
\rjpundemi= \rjmundemi + \Delta r,
$
satisfying  
$\rjmundemi = j \Delta r,\quad \rjjj= (j+1/2)\Delta r.
$
Moreover we denote by $\Delta t$ the constant time length and we set  $t_n= n \Delta t$.

As  just explained, the finite volume method is based on averaging the balance law (\ref{HBL}) over each grid cell $[ t_n, t_{n + 1}]\times I_j$, which yields
$$
\int_{I_j \times [t_n,t_{n+1}]}\mbox{div}^\omega (T(v)) \,\overline\om \, dx^0 dx^1
= \int_{I_j \times [t_n,t_{n+1}]}  S (v) \,\overline\om \, dx^0 dx^1
$$
and rearranging the terms gives
\begin{eqnarray*}
\int_{I_j \times [t_n,t_{n+1}]}  S(v) \,\overline\om \, dx^0 dx^1
= \, && \int_ {\rjmundemi}^{\rjpundemi}
 \big( T^0(t_{n+1},r) -  T^0(t_{n},r)\big)\, \overline{\om}\,dx^1
\\
&& + \int _{t_n}^{t_{n+1}}\big(\overline{\om} (\rjpundemi)T^1(t,\rjpundemi) -
\overline{\om}(\rjmundemi)T^1(t,\rjmundemi)\big)dx^0.
\end{eqnarray*}
We
recall that  $\omega$ is independent of the time variable,
so we can write 
\begin{eqnarray*}
\int_ {\rjmundemi}^{\rjpundemi}  T^0(t_{n+1},r) \, \overline{\om} dx^1
=
&&
\int_ {\rjmundemi}^{\rjpundemi} T^0(t_{n},r) \,\overline{\om} dx^1 
+
\int_{I_j \times [t_n,t_{n+1}]}  S(v)\,\overline\om \, dx^0 dx^1
\\
&& - \Bigg( \overline{\om}_{j + 1/2} \int _{t_n}^{t_{n+1}}T^1(t,\rjpundemi) dx^0-
\overline{\om}_{j - 1/2} \int _{t_n}^{t_{n+1}}T^1(t,\rjmundemi) dx^0          \Bigg).
\end{eqnarray*}
We next introduce the following  approximations of numerical flux functions and the source term
$$
\frac{1}{\Delta t}\int _{t_n}^{t_{n+1}}T^1(t,r_{j   \pm 1/2}) dx^0 \simeq   Q_{j\pm 1/2}^n,
\quad  \frac{1}{\Delta r}\int_ {\rjmundemi}^{\rjpundemi}  T^0(t_{n},r) \, \overline{\om}\, dx^1 \simeq \overline{T}_{j}^n \, \overline{\om}_j
$$
and
$$
\frac{1}{\Delta r \,\Delta t}\int_{I_j \times [t_n,t_{n+1}]}S(v) \,\overline\om dx^0 dx^1\simeq   \overline{S}_j^n\,  \overline{\om}_j.
$$

Hence, based on these  approximations,  the finite volume approximation takes the form 
$$
\overline{\om}_{j} \,  T^{n+1}_j= \overline{\om}_{j} \,  \overline{T}^n_j -
\lam \Big( \overline{\om}_{j + 1/2}\, Q_{j + 1/2}^n  -
\overline{\om}_{j - 1/2}\, Q_{j - 1/2}^n\Big) +\Delta t\,\overline{S}_{j}^n\, \overline{\om}_j,
$$
where $\lam:= \Delta t/ \Delta r. $ Since
$\overline{T}^n_{j} =\overline{T} (v_j^n) $ is invertible, after dividing the above equation by $\overline{\om}_j$, we 
arrive at the scheme
\be
\label{scheme}
v^{n+1}_j= {\overline{T}} ^{-1} \Bigg( \overline{T}(v^n_j) - \frac{\lam}{\overline{\om}_j} \Big( \overline{\om}_{ j + 1/2}\, Q_{j + 1/2}^n  -
\overline{\om}_{j - 1/2}\, Q_{j - 1/2}^n\Big) + {\Delta t\,} \overline{S}(v^n_{j})\Bigg).
\ee


\subsection{Well--balanced scheme for hyperbolic  equations on Schwarzschild spacetime}

We can think of a general solution to (\ref{HBL}) as being, locally in space and time, a small perturbation of a steady solution, that is,
in a regime where the source term balances the flux gradient within a certain moving frame. Asymptotically in time, one expects the behavior to be steady, since the geometry is static, as is the case for the models in the previous section.   
We focus attention on defining a  well balanced scheme specifically in the case 
$$
\om (r) = r( r - 2m), 
$$
for some mass $m \geq 0$, and we present the scheme in the case of the Burgers equation on Schwarzschild spacetime, that is, (\ref{RBE}).
The discrete version of (\ref{RBE}) now reads 
\be
\label{DS1}
\overline{T}_j^{n+1}= \overline{T}_j^n - \frac{\Delta t}{\Delta r}(\om_{j+1/2}Q_{j+1/2} - \om_{j-1/2}Q_{j-1/2})+ \Delta t \, \overline{S}_j, 
\ee
where $\Delta r$ still denotes the mesh size $\Delta r= \rjpundemi-\rjmundemi$, 	
but we now take into account the mass parameter in order to define the mesh points, that is, 
$$
\rjmundemi = 2m +j \Delta r, \quad \rjjj= 2m +(j+1/2)\Delta r, \quad \rjpundemi = 2m +(j+1) \Delta r
$$
and the averaged weights are 
$
\om_{j\pm1/2}= r_{j\pm1/2}(r_{j\pm1/2}- 2m).
$
Moreover, $\overline{T}_j^n$ and $\overline{S}_j$ are approximations corresponding to 
$$
\overline{T}_j^n=\frac{1}{\DeltaVt} \int_{\rjmundemi}^{\rjpundemi} T \,dV_{\widetilde{g}},
\qquad 
\overline{S}_j= \frac{1}{\DeltaVt}  \int_{\rjmundemi}^{\rjpundemi}  S \, dV_{\widetilde{g}}, 
$$
in which $\widetilde{g}$ is the induced metric on spacelike slices of the foliation, with $g = -c^2\,\bigl( 1 - {2m \over r} \bigr) \, dt^2 + \widetilde{g}$ and $\dVt$ is the natural volume form induced on a spacelike slice, so 
$$
\dVt := r^2\, \Big(\frac{r}{r-2m} \Big)^{1/2} \,dr, \qquad
\DeltaVt :=r^2\, \Big( \frac{r}{r-2m} \Big)^{1/2} \, \Delta r.
$$

It is convenient to introduce $\Omega :=\om\, v^2/2$. To treat the source term, we integrate the equation (\ref{SS3})  over a spacetime cell, and obtain 
$$
\int_{r-1/2+}^{r+1/2-} \sqrt{\widetilde{g}}\, \Omega dr
= 
\int_{r-1/2+}^{r+1/2-} (rv^2- m c^2) \, \sqrt{\widetilde{g}}  \,dr = \DeltaVt \overline S_j, 
$$
which gives
$
\overline S_j= \frac{1}{\DeltaVt} \int_{r-1/2+}^{r+1/2-} \sqrt{\widetilde{g}}\, \Omega dr.
$
By substitution and integration by parts, we arrive at the expression of the source term: 
$$
\overline S_j
= \frac{1}{\DeltaVt} \Big( r^{7/2}(r-2m)^{1/2}v^2/2\Big)\Big|_{r-1/2+}^{r+1/2-}
-{1\over \DeltaVt}\int_{r-1/2+}^{r+1/2-} \frac{r(r-2m)}{2}\, v^2\, \Bigg({2\over r} - {m \over r(r-2m)}\Bigg) \, dV_{\widetilde g}. 
$$
In the numerical experiements, we will use Gaussian quadrature to evaluate this integral.

Furthermore, the left-- and right--hand numerical flux terms are given by 
$$
Q_{j+1/2}=Q\big(\overline{T}_{j+1/2-}, \overline{T}_{j+1/2+}\big), \qquad
Q_{j-1/2}=Q\big(\overline{T}_{j-1/2-}, \overline{T}_{j-1/2+}\big),
$$
where
$\overline{T}_{j+1/2\pm} = r^2_{j+1/2}\overline{v}_{j+1/2\pm}$ 
with reconstructed velocities 
\begin{eqnarray*}
&\overline{v}_{j+1/2+}= \pm \sqrt{c^2-K_1^2\, \Big(1-\frac{2m}{\rjpundemi}\Big)}, \qquad
\overline{v}_{j+1/2-}= \pm \sqrt{c^2-K_2^2\, \Big(1-\frac{2m}{\rjpundemi}\Big)}, 
\\
&\overline{v}_{j-1/2+}= \pm \sqrt{c^2-K_2^2\, \Big(1-\frac{2m}{\rjmundemi}\Big)},\qquad
\overline{v}_{j-1/2-}= \pm \sqrt{c^2-K_3^2\, \Big(1-\frac{2m}{\rjmundemi}\Big)},
\end{eqnarray*}
and 
\begin{eqnarray*}
K_1^2= \frac{c^2-\overline{v}^2_{j+1}}{1- \frac{2m}{\rjjj +\Delta r}}, \qquad
K_2^2= \frac{c^2-\overline{v}^2_{j}}{1- \frac{2m}{\rjjj}},\qquad
K_3^2= \frac{c^2-\overline{v}^2_{j-1}}{1- \frac{2m}{\rjjj -\Delta r}}.
\end{eqnarray*}
This reconstruction is motivated by the steady solutions derived in (\ref{SS2}) which allow us to compute intermediate states at which the numerical flux is evaluated. Importantly, in the special case that $K_1, K_2, K_3$ coincide, the proposed scheme recovers an {\sl exact steady solution.} 


\subsection{Second-order accurate version} 

To increase the order of accuracy of the above scheme, we rely on the technique introduced by Nessyahu and Tadmor central scheme \cite{NessyahuTadmor}. Our second--order version is thus a predictor-corrector scheme, defined by 
\begin{eqnarray*}
&v^{n+1/2}_{j}=v^{n}_{j}-{1\over 2}\, \lambda \,\del_0T_{n,j}^1, 
\\
&v^{n+1}_{j}={1\over 2}(v^{n}_{j+1}+v^{n}_{j-1})+{1\over 4}(v'_{j-1}-v'_{j+1})
-\frac{\lambda}{2}\bigg({T^1}(v^{n+1/2}_{j+1})-{T^1}(v^{n+1/2}_{j-1})\bigg),
\end{eqnarray*}
where $\lambda :=\Delta   t/\Delta  r$ denotes the mesh ratio and $\frac{1}{\Delta x} \,\del_0T_{n,j}^1$ is the approximate first--order
derivative of the flux $T_{n,j}^1= T^1( v_{n,j})$. The time step must satisfy the stability condition
$
\lambda \,\max_j\,|\del_0T_{n,j}^1|\leq\frac{1}{2},
$
where $|\del_0T_{n,j}^1|$ denotes the spectral radius of the Jacobian matrix. In order to get rid of spurious oscillations, a slope limiter must be applied to the values $v_j'$ and $\del_0T_{n,j}^1$ and, specifically, 
\begin{eqnarray*}
v_{j}' &=& \mbox{Minmod}\,(v _{j+1} - v _{j},\,v _{j} -v_{j-1}),
\\
\del_0T_{n,j}^1 &=& \mbox{Minmod}\, \big( T^1_{j+1} - T^1_j ,\, T^1_j -T^1_{j-1} \big),
\\
\mbox{Minmod}(a,b) &=&  
{\rm sgn}(a) \, \min(|a|, |b|) \, \, {\rm if } \, \,  ab>0; 
\qquad  
 0 \, \, {\rm if }   \, \, ab \leq 0.
\end{eqnarray*}


\section{Numerical experiments}

\subsection{Comparison between several schemes and models}

Numerical experiments are now presented for both models derived in this paper. After normalization ($c=1/\eps=1$), we thus consider the  
\mbox{geometric Burgers equation (I):} 
\begin{eqnarray}
\label{ss5}
 \del_t (r^2 \, w) + \del_r \bigg( r(r-2m) \, \big( -1 + \sqrt{ 1+ w^2} \big) \bigg) = 0, 
\end{eqnarray}
and the \mbox{geometric Burgers equation (II):} 
\begin{eqnarray}
\label{RB5}
 \del_t (r^2\,v) + \del_r\bigg(r(r-2m)\, \frac{v^2}{2}\bigg) = rv^2 - m,
\end{eqnarray}
with obvious notation. The first equation has a conservative form, while the second one is nonconservative. 
We work within the interval $r\in (2m, R)$, where $R$ is a fixed upper bound for the spatial variable. We impose no-influx boundary conditions at $r=2m$ and $r=R$ by solving a Riemann problem between the boundary data determined by the initial data and the current numerical value at the boundary. In other words, we use the Godunov flux at the boundary.  
For all numerical experiments, we use $2m=0.1$ and $R=1.0$ and the value of CFL number is taken to be $0.9$.
We begin by illustrating the importance of the well balanced property by comparing
together the numerical solutions obtained with the proposed scheme and with a more naive discretization of
the right--hand side of the geometric Burgers equation (II).
For the model (I), we used $\Delta r=0.02$ and the endpoint velocity $v(R) =0.03$. The numerical solution is plotted on Figure~5.1.1.
In this Figure, in (b), (c), we plot the numerical solutions given by the three schemes. The first one
      is the first--order Lax--Friedrichs scheme (which is indicated by $+$); the second one is the second--order 
     version of Lax--Friedrichs scheme (which is indicated by line); the third one is the well balanced second--order version of Lax--Friedrichs scheme (which is indicated
      by $-$). In Figures (d), (e), we compare the numerical flux $r(r-2m) \, f_\eps(w)$ (cf.~(\ref{ss5})), 
      by the three schemes.
The proposed well balanced scheme appears to be efficient and robust and, in particular, this scheme preserves steady solutions.
Next, for the model (II), we used $\Delta r=0.005$ and the velocity $v(R) =0.3$, which led us to Figure~5.1.2.
In Figures~(b) and (c), we compare the numerical solutions based on the three schemes.
The first one is a standard first--order Lax--Friedrichs scheme (which is indicated by +); the second
one is a standard second--order Lax--Friedrichs scheme (and is plotted with dots); the
third one is a well balanced second--order Lax--Friedrichs scheme (and is plotted with $-$). In Figures (d) and (e), we compare the numerical values of $K^2$ (defined in (\ref{SS2})) based on the three
schemes. Again, the scheme is found to be efficient and preserve steady solutions.

\begin{figure}[t!]
\centering 
\subfloat[]{
{
    \begin{minipage}[b]{0.48\linewidth}
        \centering \includegraphics[width=0.98\linewidth]{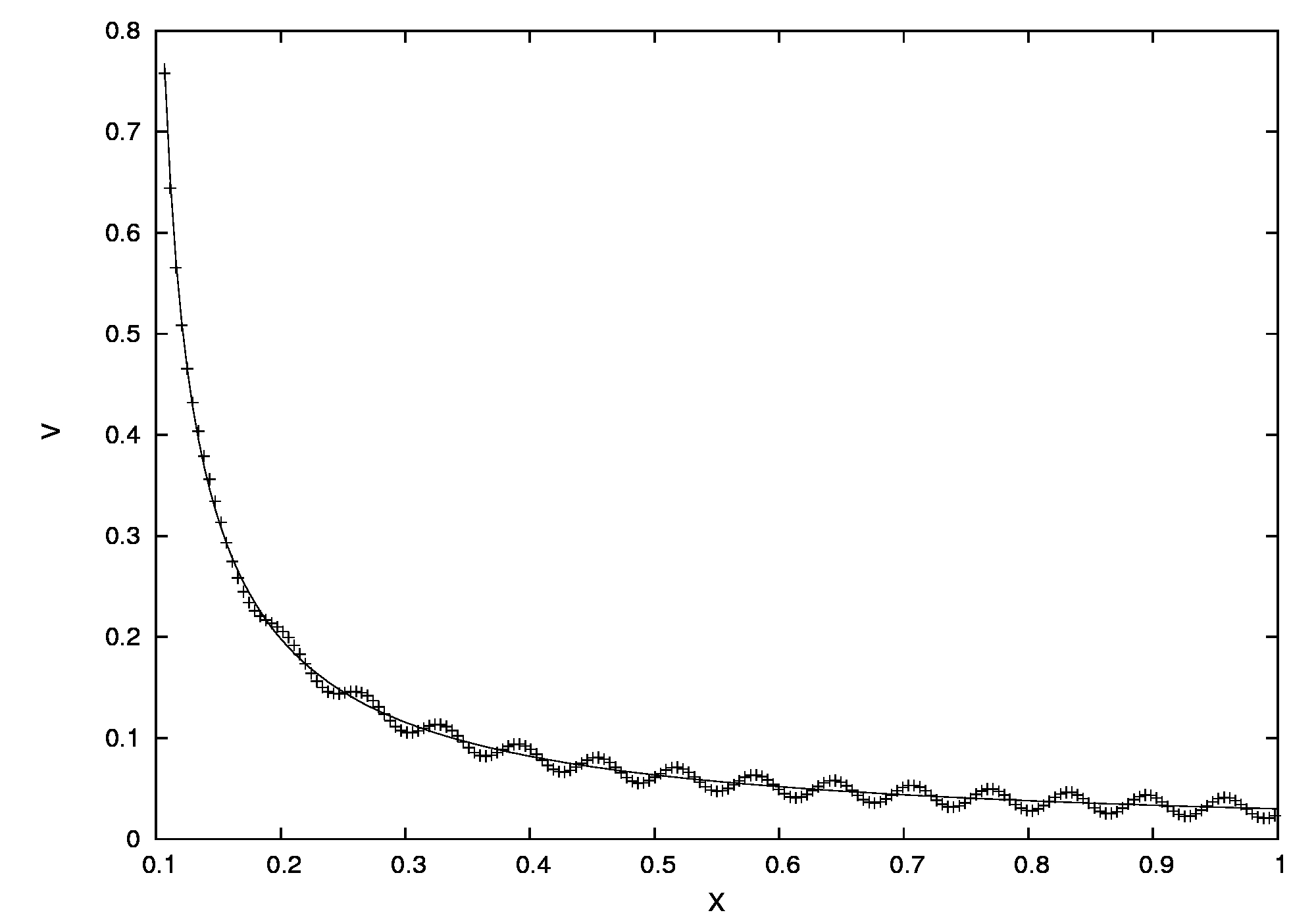}
    \end{minipage}}}
\\
\subfloat[]{
    \begin{minipage}[b]{0.48\linewidth}
        \centering \includegraphics[width=0.98\linewidth]{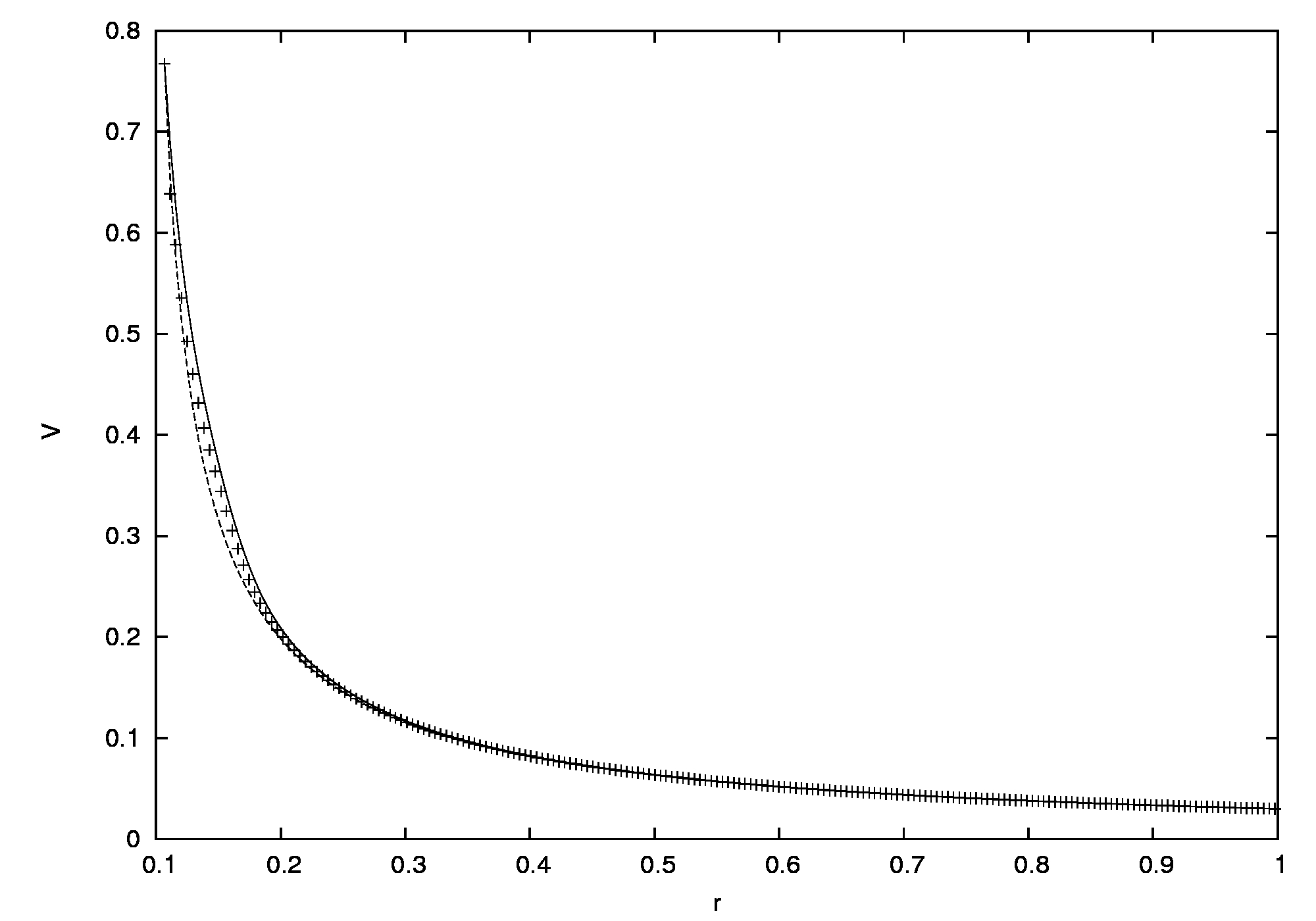}
    \end{minipage}
} \hfill \subfloat[]{
    \begin{minipage}[b]{0.48\linewidth}
        \centering \includegraphics[width=0.98\linewidth]{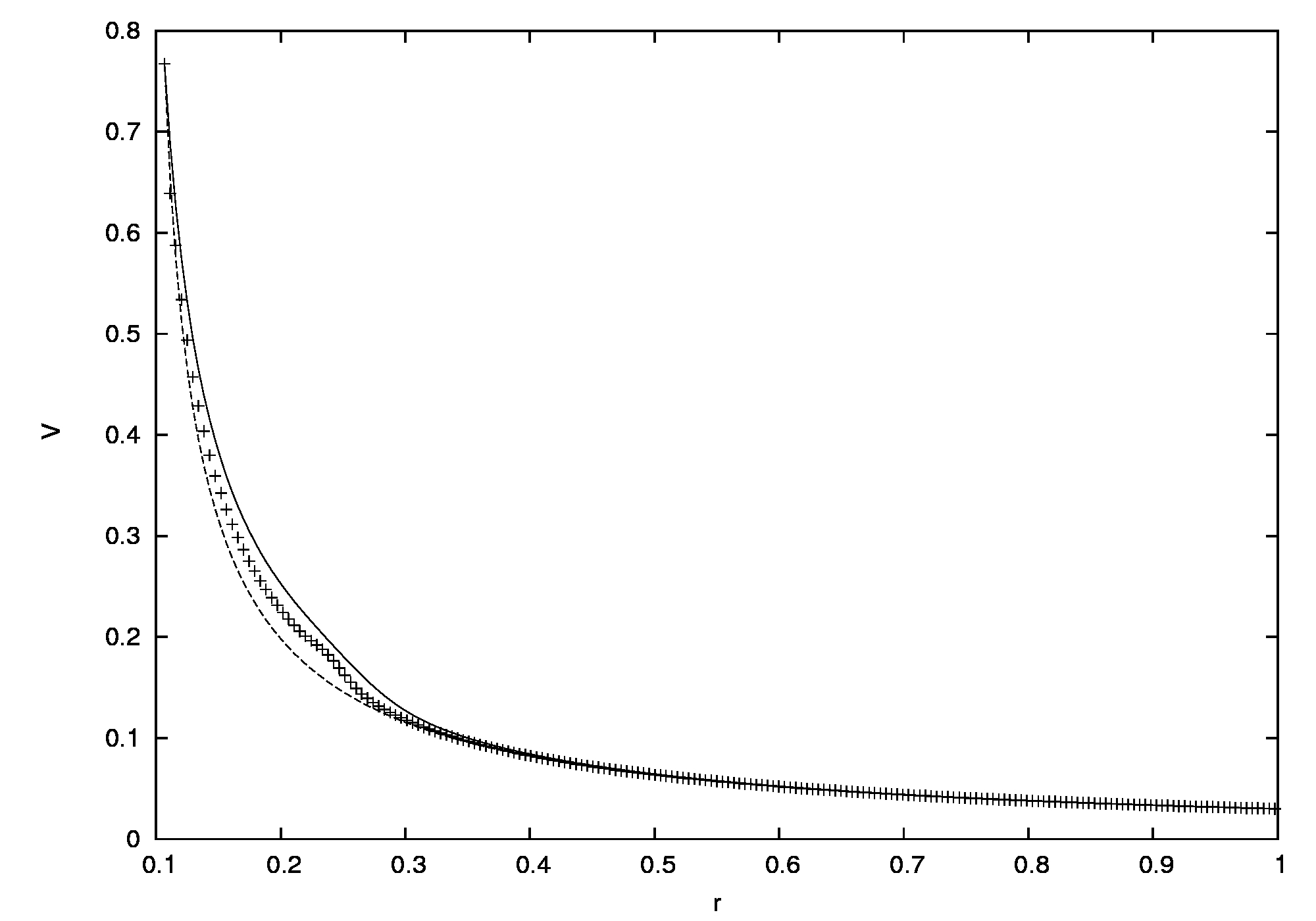}
    \end{minipage}
} \hfill \subfloat[]{
    \begin{minipage}[b]{0.48\linewidth}
        \centering \includegraphics[width=0.98\linewidth]{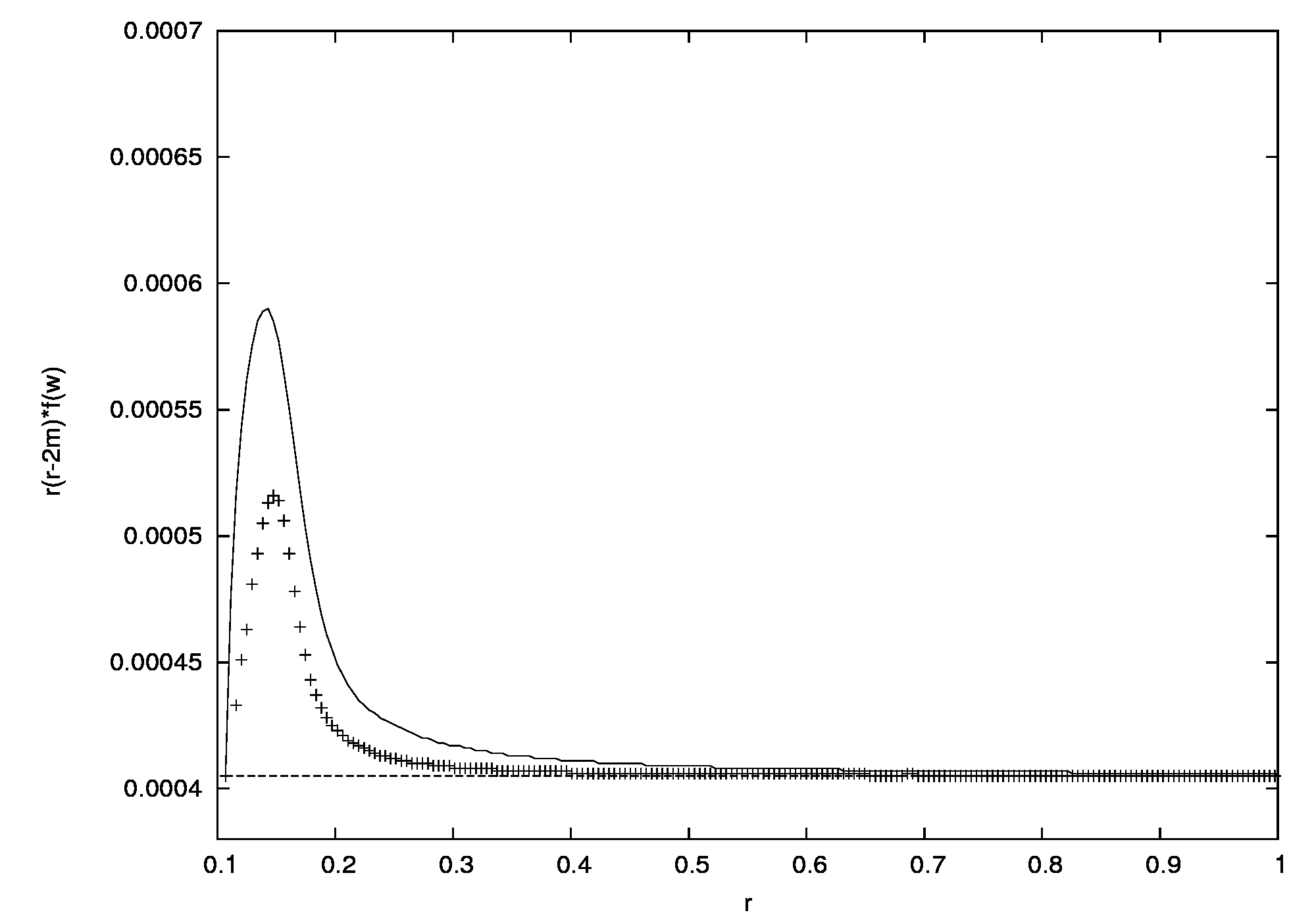}
    \end{minipage}
} \hfill \subfloat[]{
    \begin{minipage}[b]{0.48\linewidth}
        \centering \includegraphics[width=0.98\linewidth]{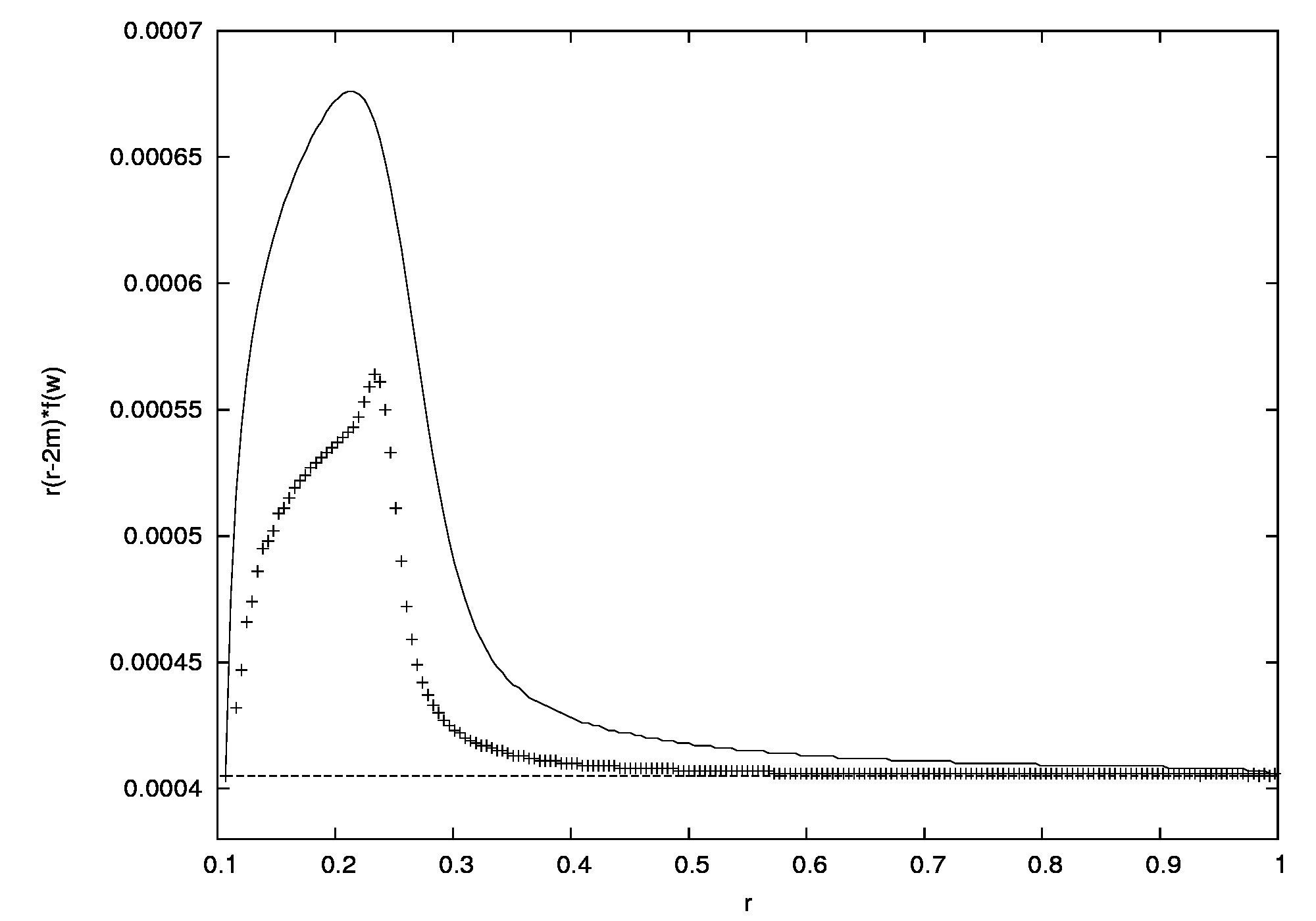}
    \end{minipage} }\begin{center}
                    {Figure 5.1.1. The numerical solutions given by the three schemes for the first model (I).}
                   \end{center}
\end{figure}

 \begin{figure}[t!]
\centering 
\subfloat[]{
{
    \begin{minipage}[b]{0.48\linewidth}
        \centering \includegraphics[width=0.98\linewidth]{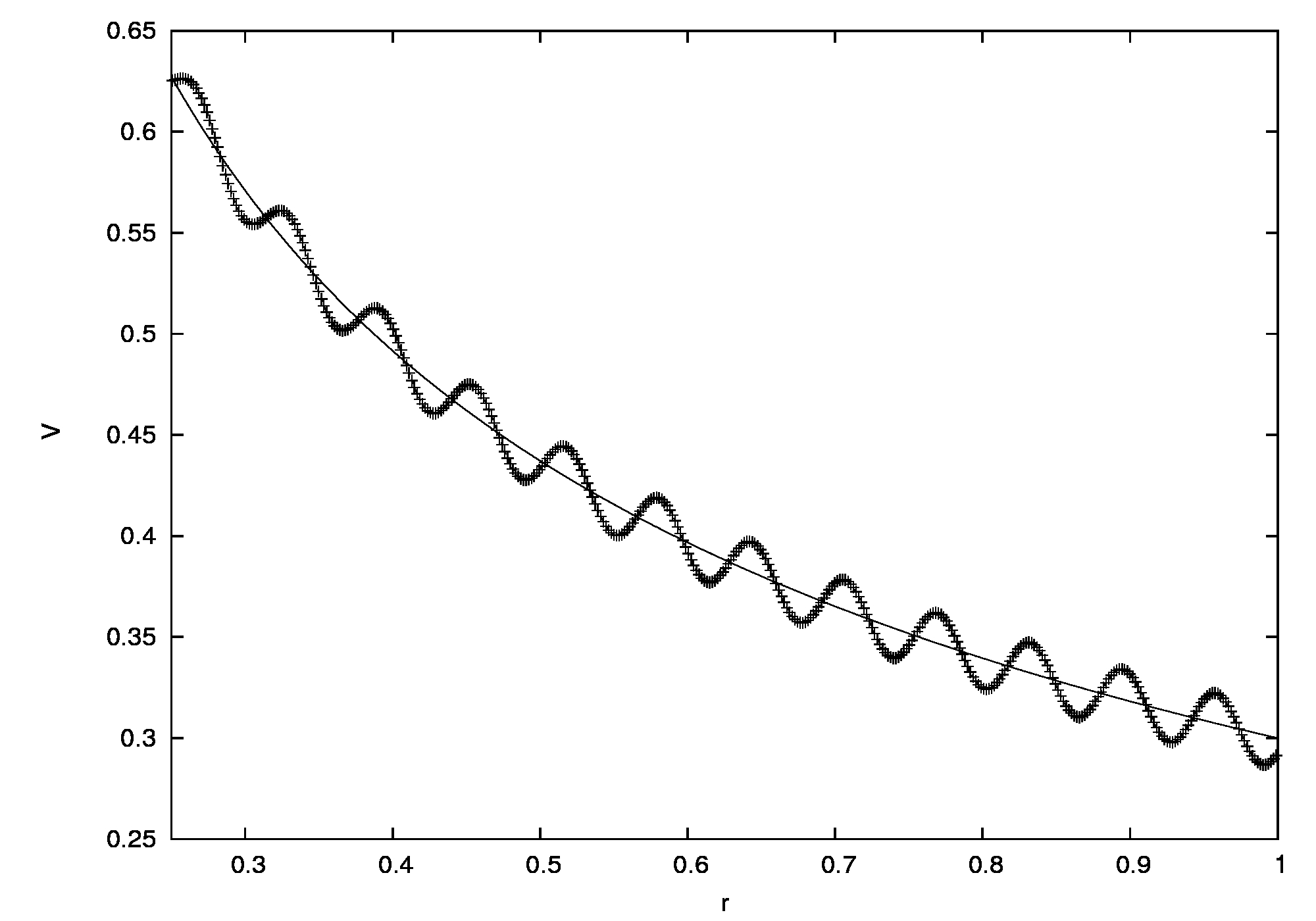}
    \end{minipage}}}
\\
\subfloat[]{
    \begin{minipage}[b]{0.48\linewidth}
        \centering \includegraphics[width=0.98\linewidth]{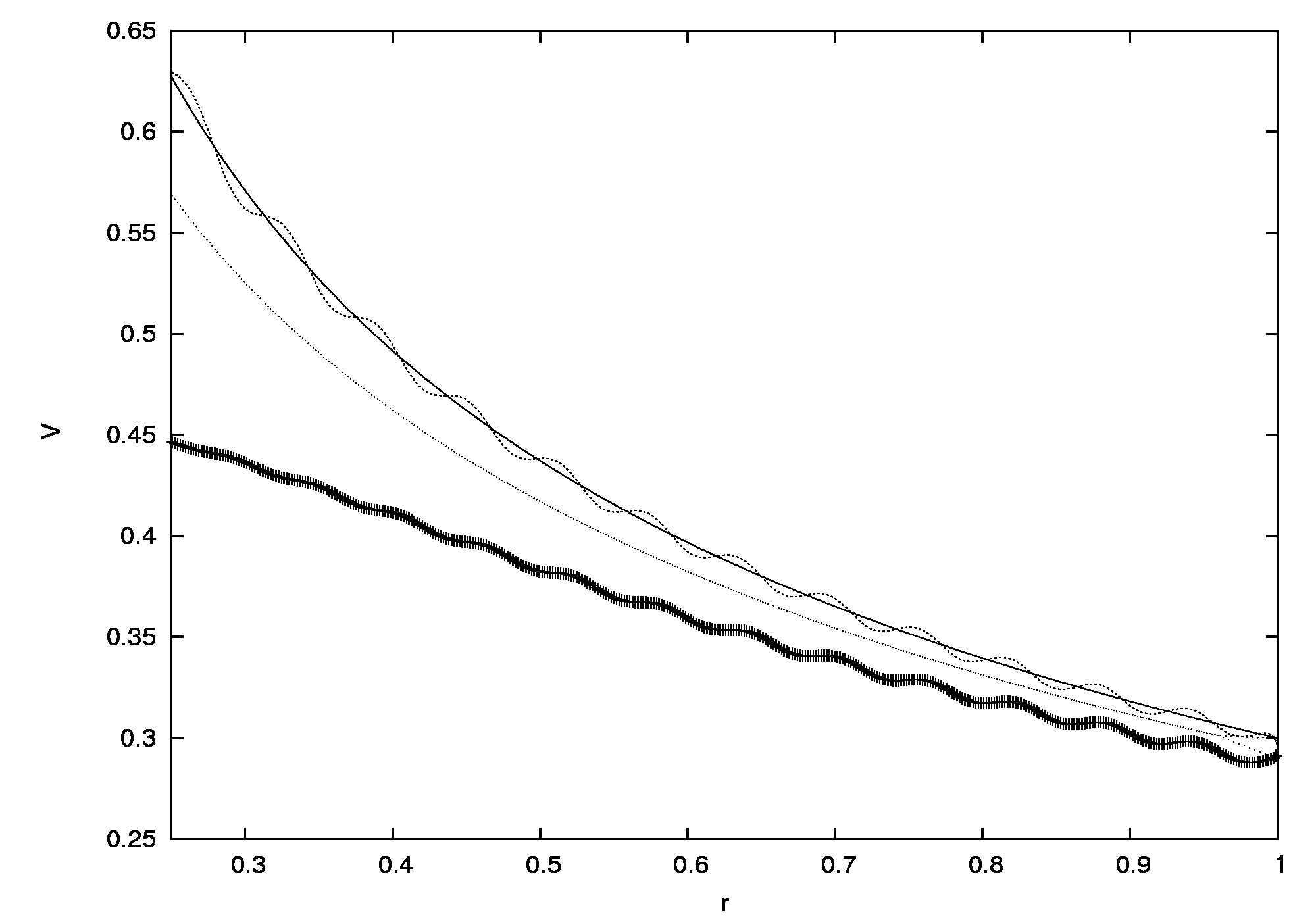}
    \end{minipage}
} \hfill \subfloat[]{
    \begin{minipage}[b]{0.48\linewidth}
        \centering \includegraphics[width=0.98\linewidth]{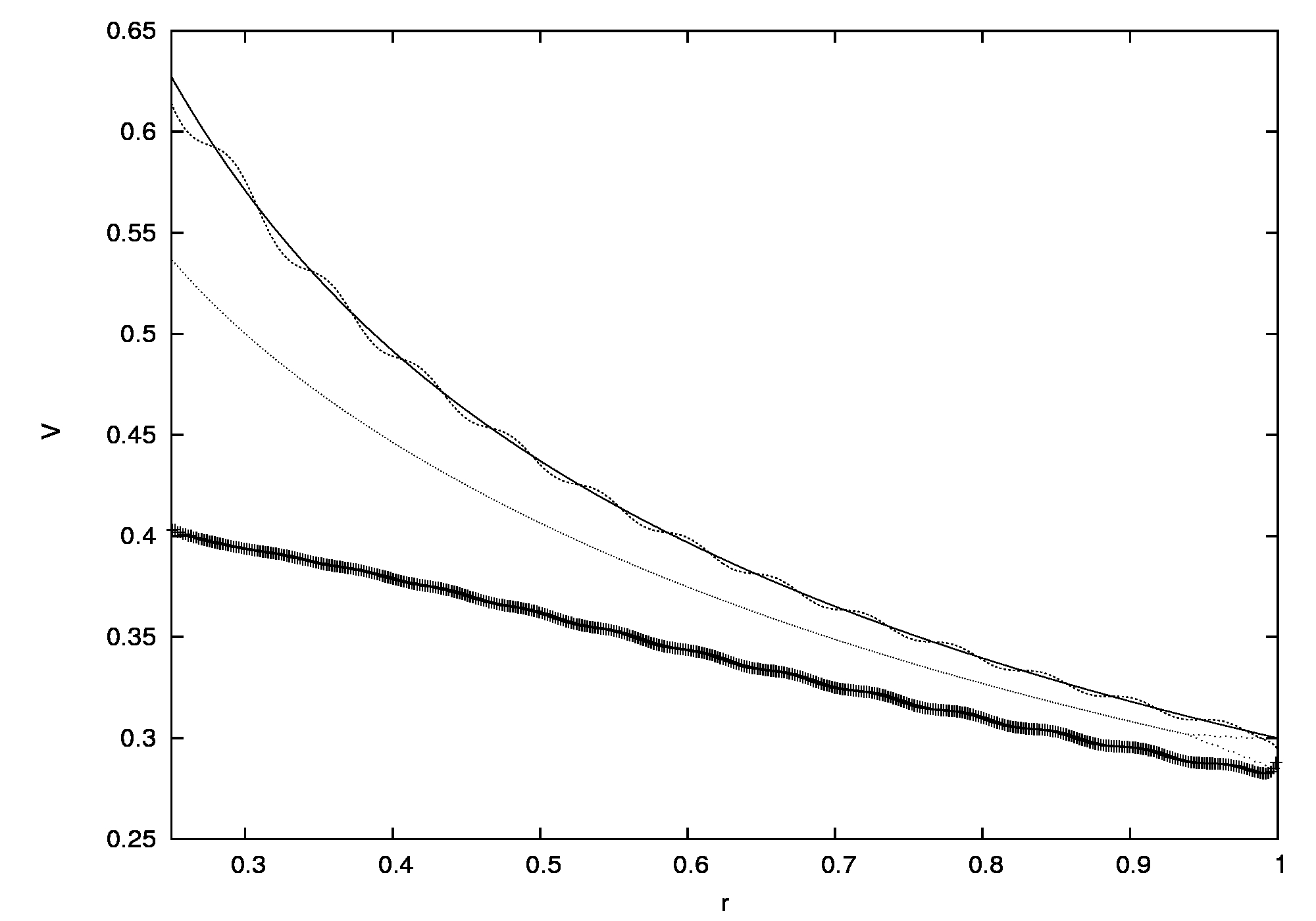}
    \end{minipage}

} \hfill \subfloat[]{
    \begin{minipage}[b]{0.48\linewidth}
        \centering \includegraphics[width=0.98\linewidth]{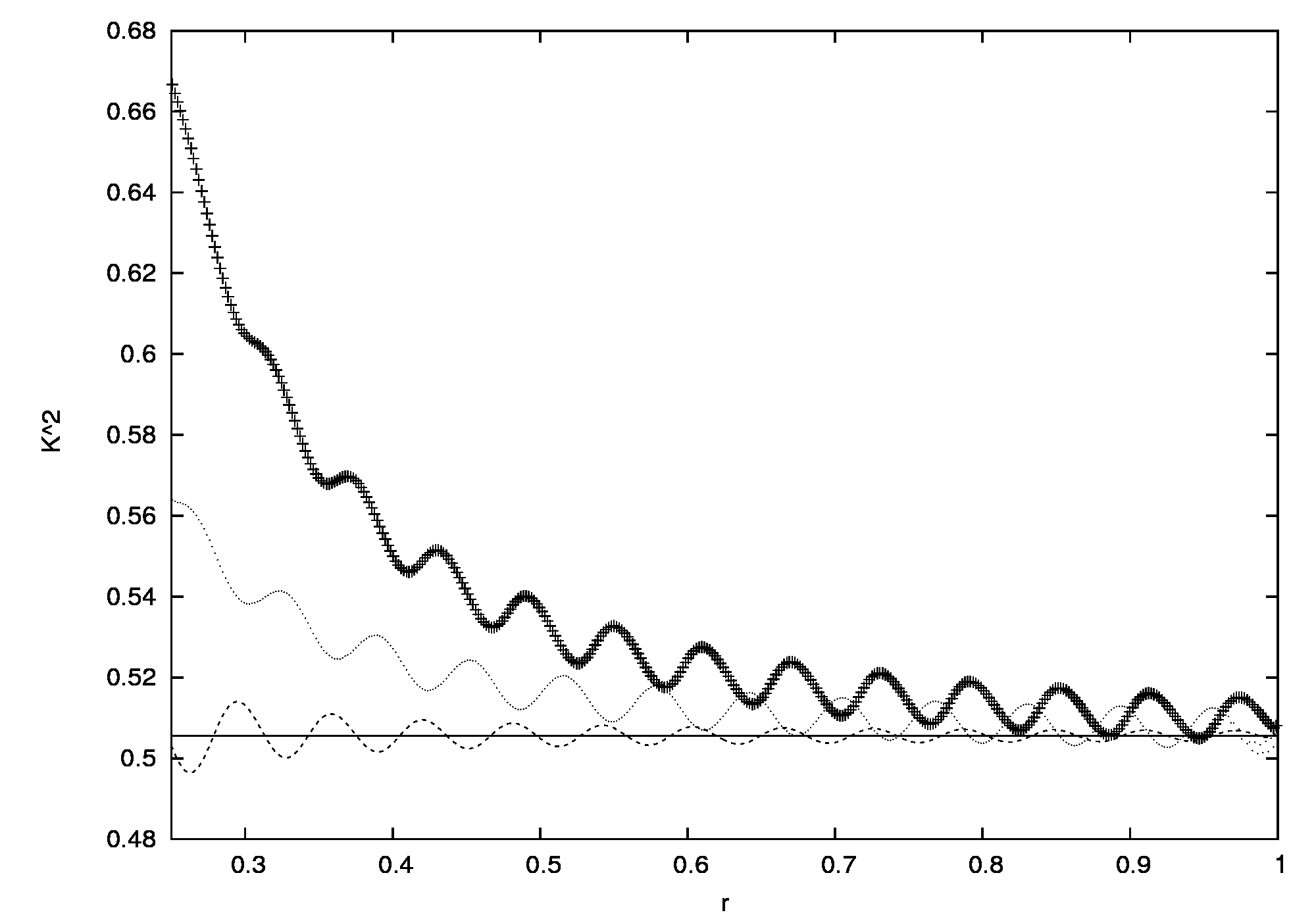}
    \end{minipage}
} \hfill \subfloat[]{
    \begin{minipage}[b]{0.48\linewidth}
        \centering \includegraphics[width=0.98\linewidth]{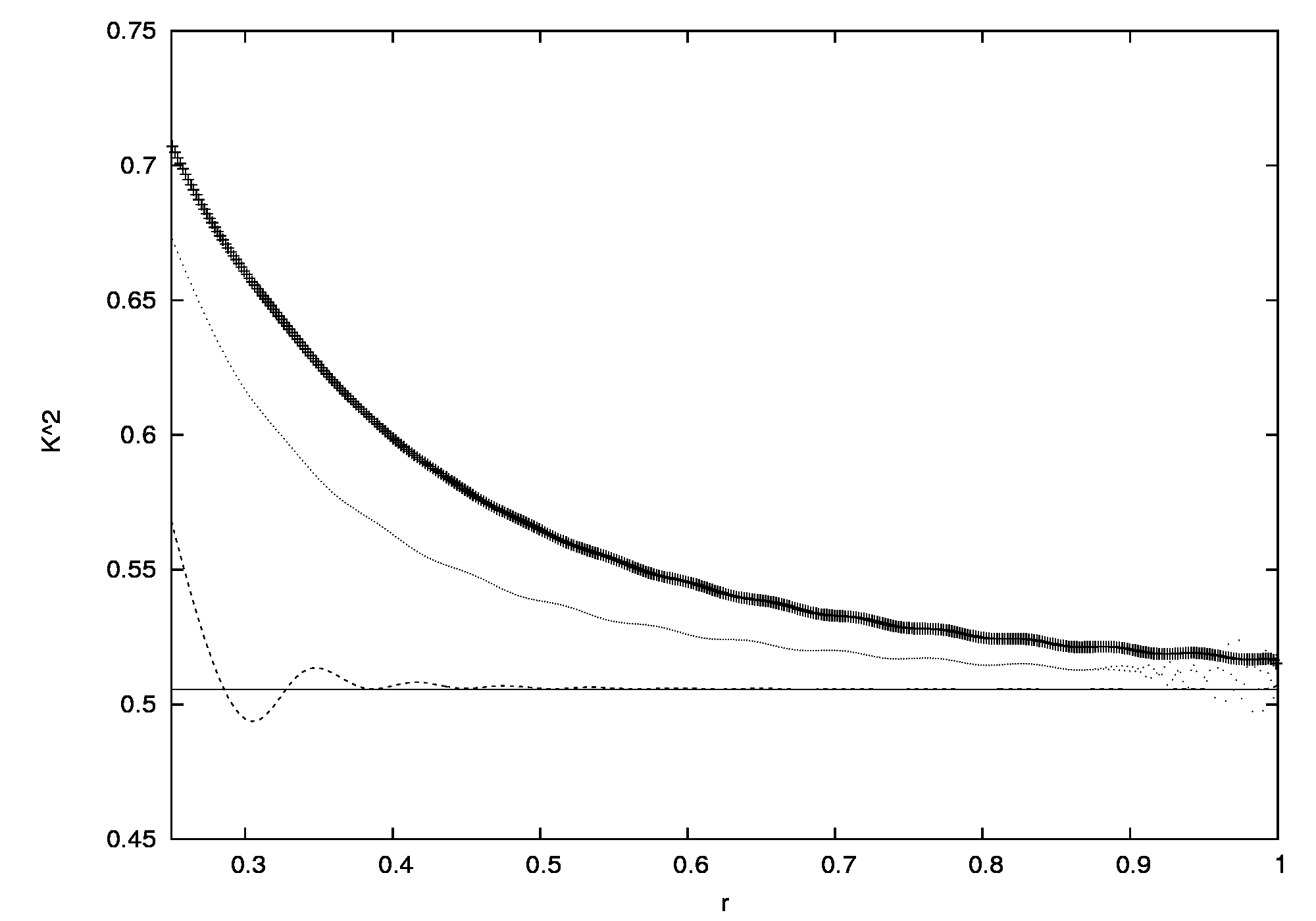}
    \end{minipage} }\begin{center}
                    {Figure 5.1.2. The numerical solutions given by the three schemes for the second model (II).}
                   \end{center}
\end{figure}


\subsection{Comparison of several schemes for a single shock}

\begin{figure}[t!]
\centering  
\subfloat[]{
    \begin{minipage}[b]{0.48\linewidth}
        \centering \includegraphics[height=120pt, width=0.98\linewidth]{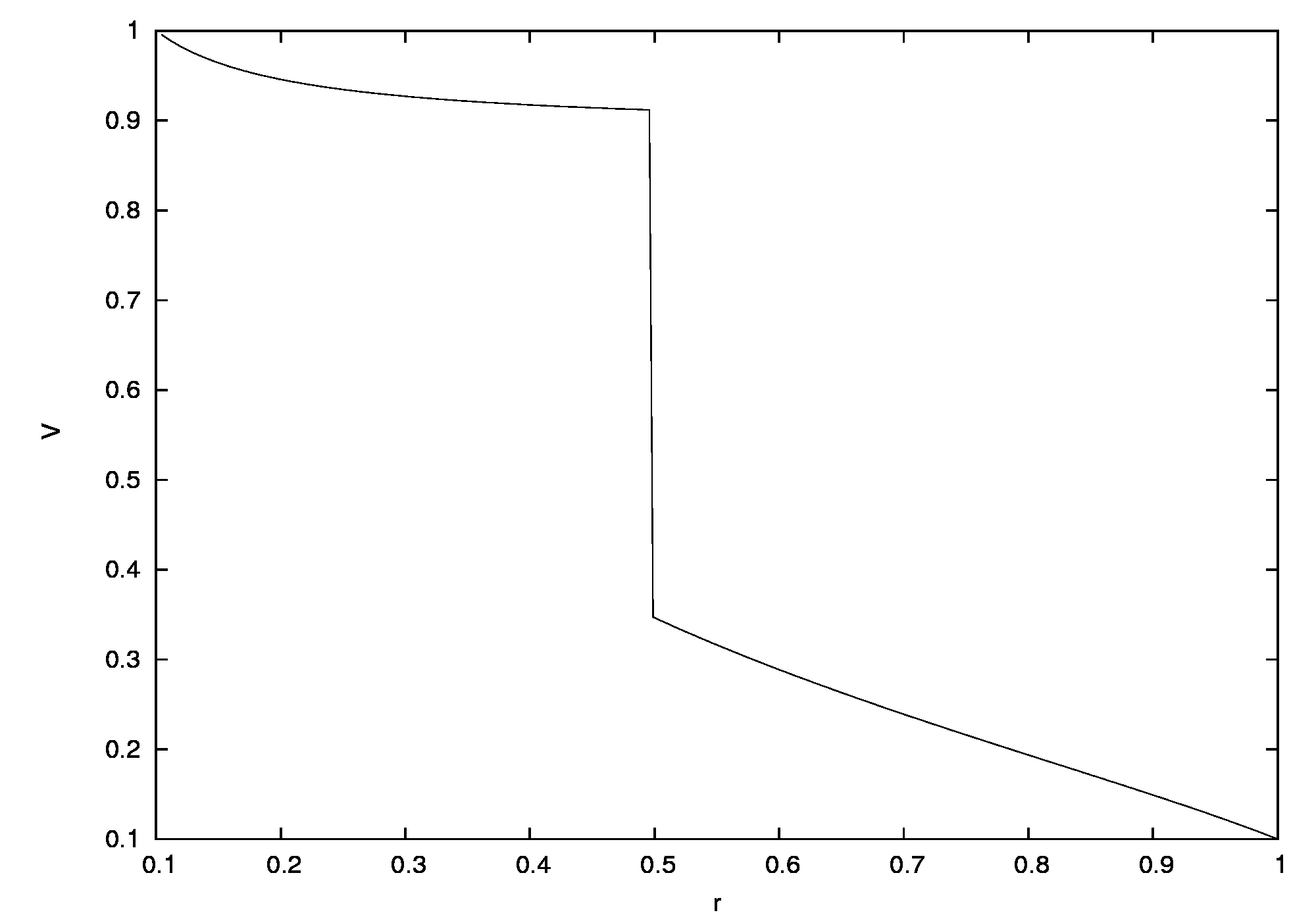}
    \end{minipage}}
\subfloat[]{
    \begin{minipage}[b]{0.48\linewidth}
        \centering \includegraphics[height=120pt, width=0.98\linewidth]{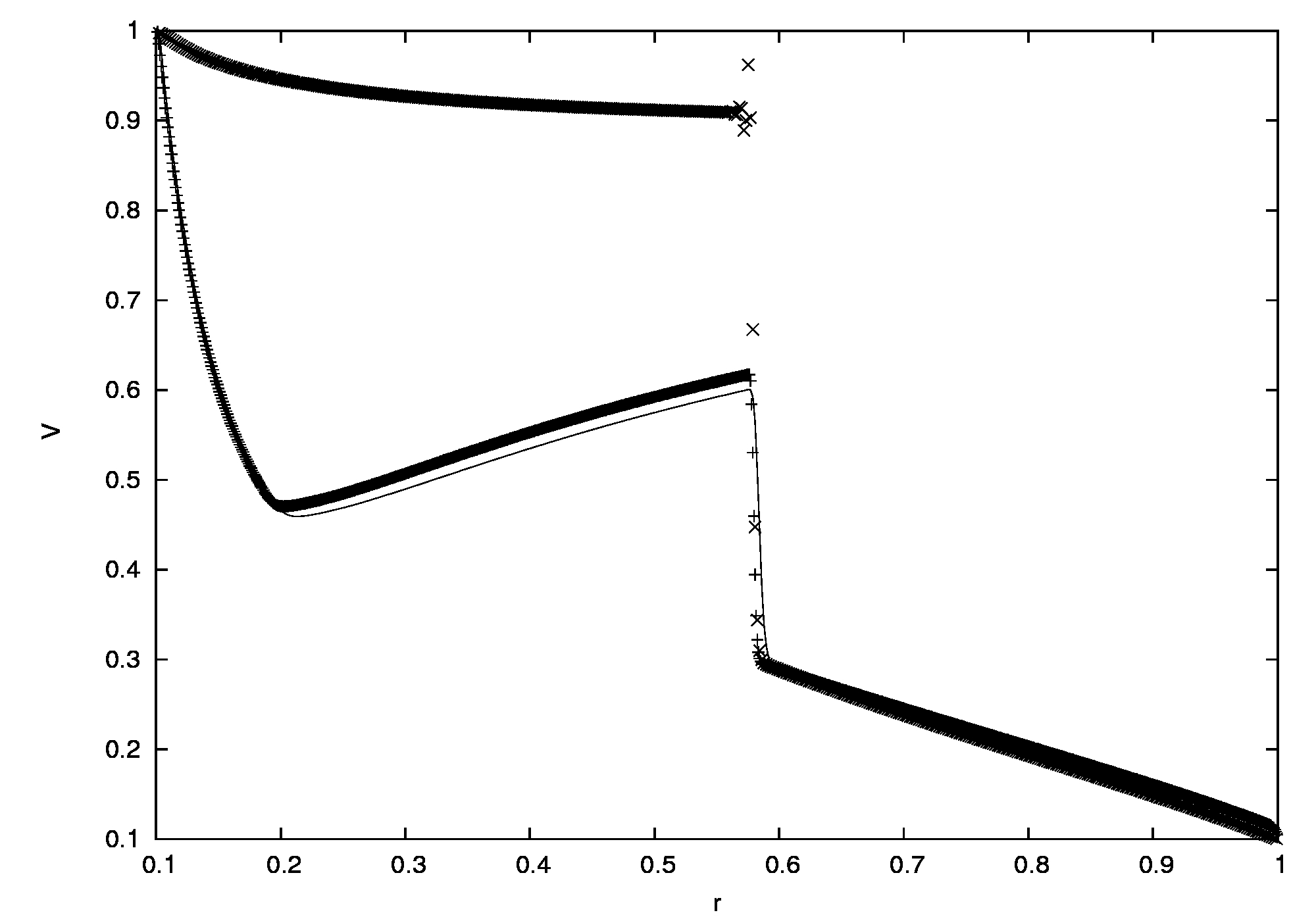}
    \end{minipage}
} \hfill \subfloat[]{
    \begin{minipage}[b]{0.48\linewidth}
        \centering \includegraphics[height=120pt, width=0.98\linewidth]{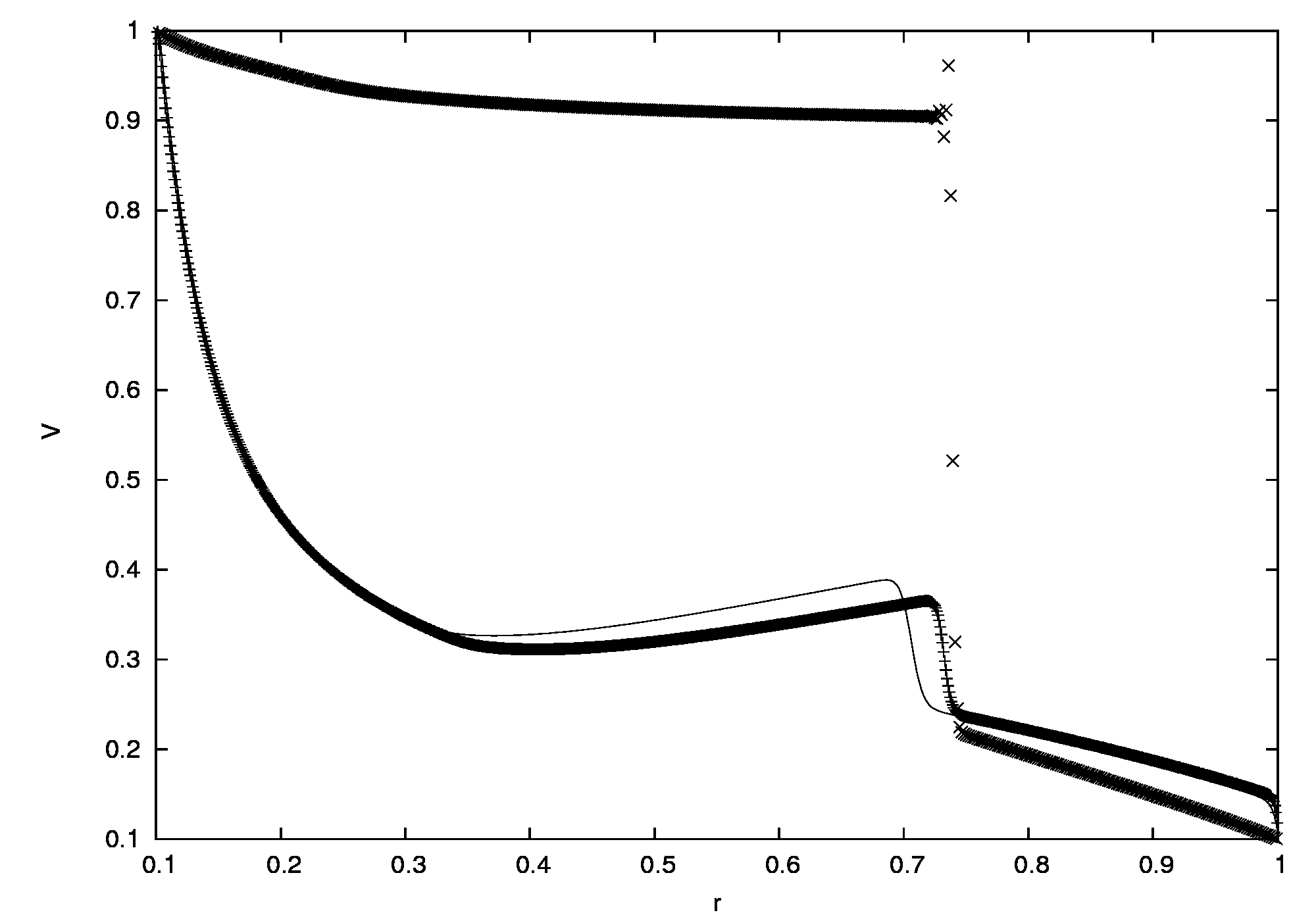}
    \end{minipage}
} \hfill \subfloat[]{
    \begin{minipage}[b]{0.48\linewidth}
        \centering \includegraphics[height=120pt, width=0.98\linewidth]{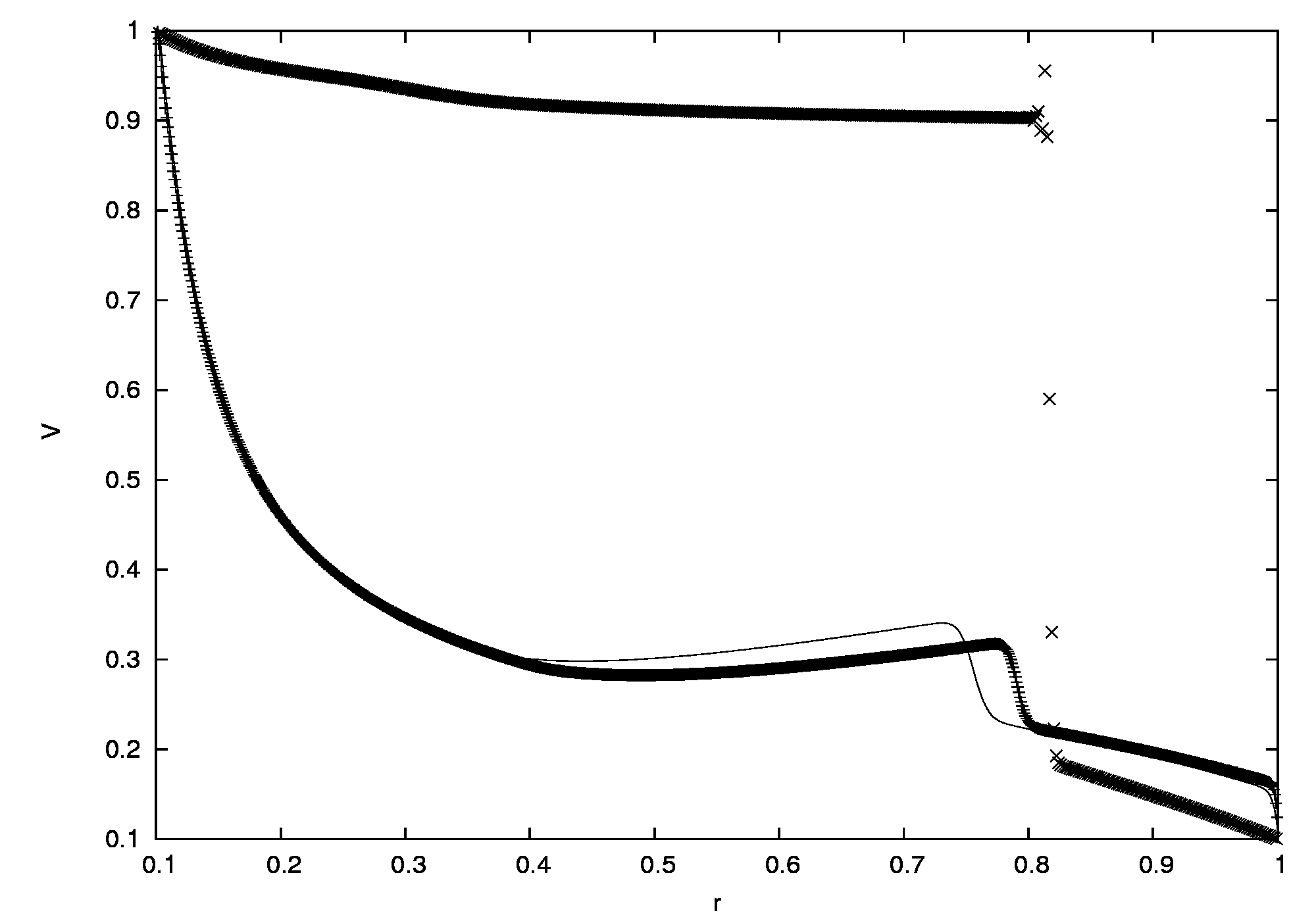}
    \end{minipage}
}\begin{center}
                    {Figure 5.2. The numerical solutions given by the three schemes for the second model (II).}
                   \end{center}
     
\end{figure}

We study next weak solutions containing a single shock separating two steady solutions at the location $r=0.5$: one steady solution corresponds to the data $v_1(R)=0.1$ while the other steady solution corresponds to $v_2(R)=0.9$. The initial data are constructed so that the initial discontinuity propagates as a single shock (without decomposing itself during the evolution), which interacts with the background geometry.  We consider the second model (\ref{RB5}) with $\Delta r=0.05$ . The numerical solution is presented  in Figure~5.2.
We compare the propagation of a single shock by three schemes. The first one is a standard first--order Lax--Friedrichs scheme (plotted with a continuous line), and the second one is a standard second--order Lax--Friedrichs scheme (ploted with the symbol $+$). Finally, the third scheme is the proposed well balanced second--order Lax--Friedrichs scheme (plotted with the symbol $\times$). One observes the efficiency and robustness of the well balanced scheme; the numerical solution based on this scheme is monotone and stable under perturbations.


\subsection{Application to late-time asymptotics -- perturbed steady solutions}

\begin{figure}[t!]
\centering \subfloat[]{
    \begin{minipage}[b]{0.48\linewidth}
        \centering \includegraphics[width=0.98\linewidth]{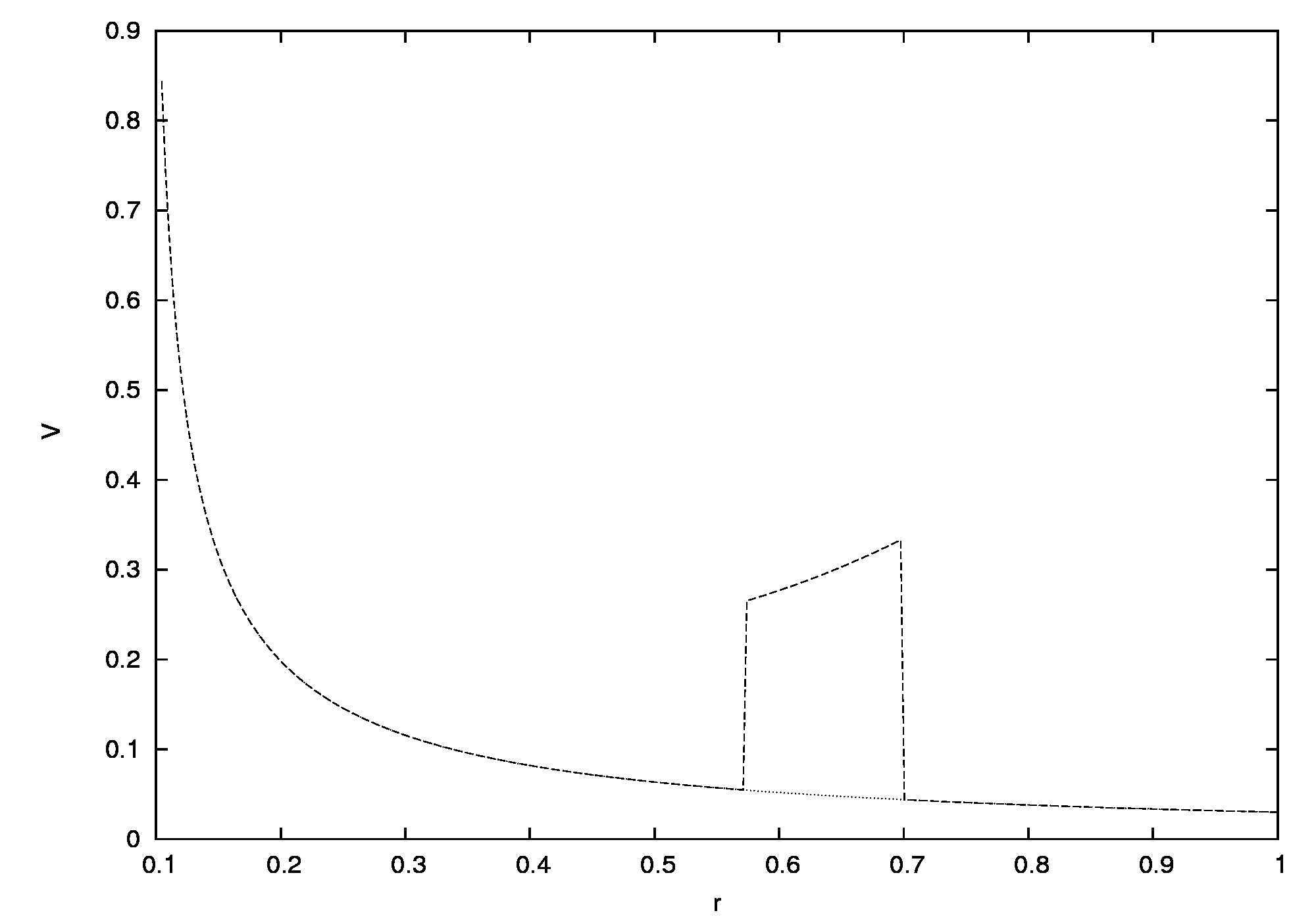}
    \end{minipage}}
\subfloat[]{
    \begin{minipage}[b]{0.48\linewidth}
        \centering \includegraphics[width=0.98\linewidth]{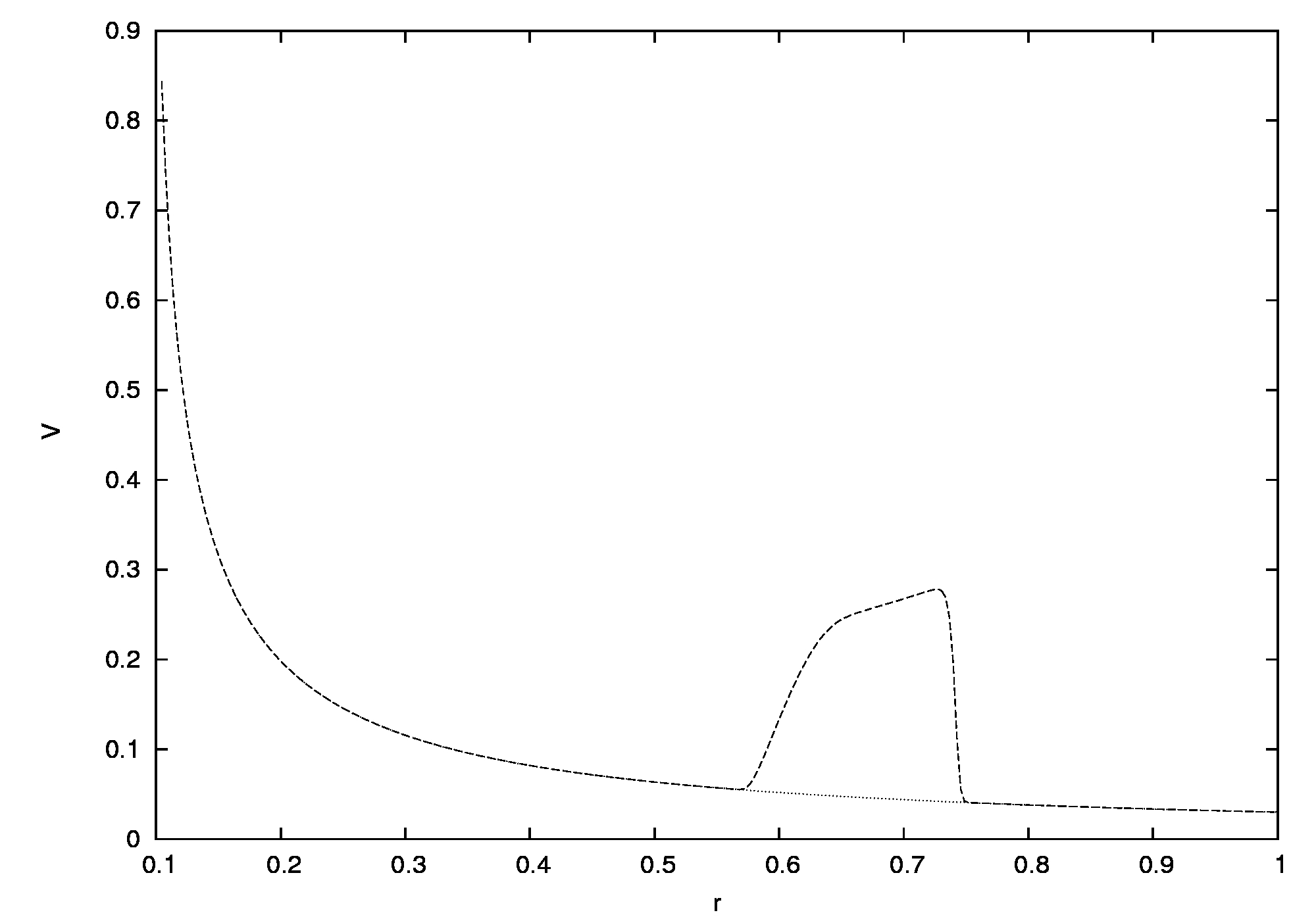}
    \end{minipage}
} \hfill \subfloat[]{
    \begin{minipage}[b]{0.48\linewidth}
        \centering \includegraphics[width=0.98\linewidth]{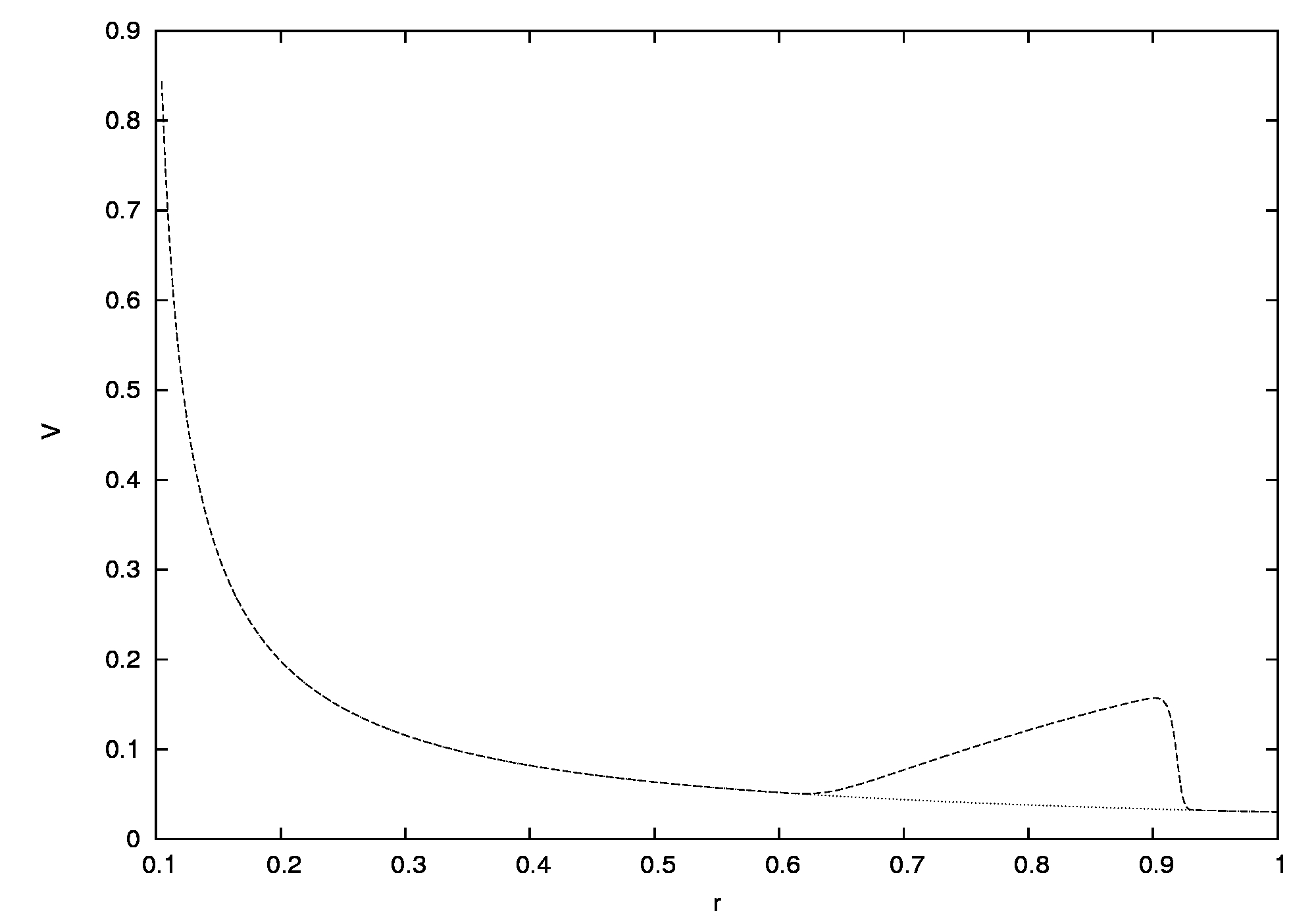}
    \end{minipage}
}
\subfloat[]{
    \begin{minipage}[b]{0.48\linewidth}
        \centering \includegraphics[width=0.98\linewidth]{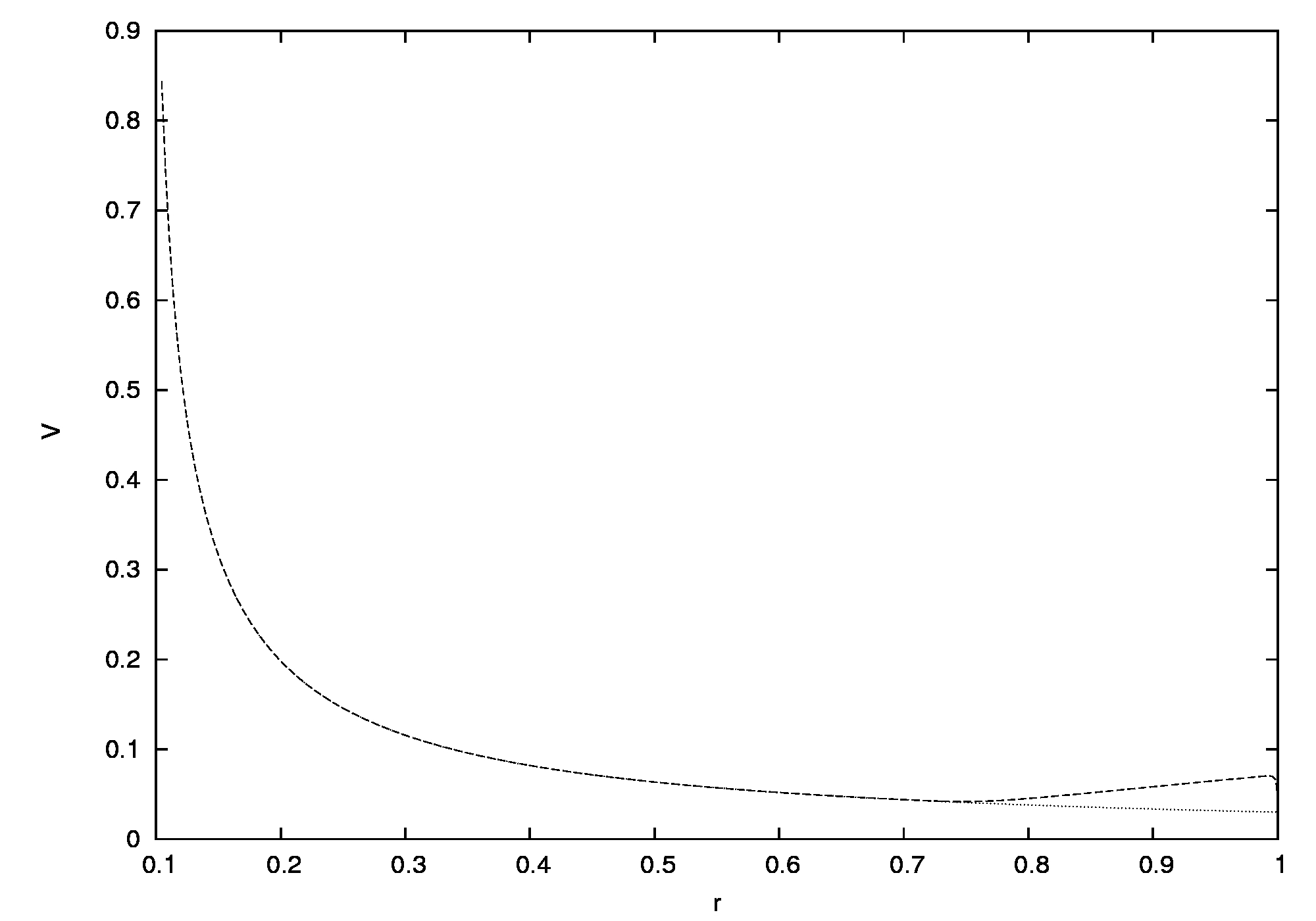}
    \end{minipage}
}\begin{center}
                    {Figure 5.3.1. Model (I). The asymptotic behaviors of long time solutions  for equation (\ref{ss5}) using the well balanced second--order scheme. }
                   \end{center}
\end{figure}

We study here the nonlinear stability of a steady solution superimposed with a compactly supported
initial perturbation.  In Figure~5.3.1, we still use (\ref{ss5}) and take each subinterval $\Delta r= 0.03$ and a velocity at the right-hand point $v(R)=0.03$. We observe that the initial perturbation evolves toward a so-called $N$--wave, which is 
typical of nonlinear hyperbolic systems, and that this perturbation moves in the right direction away from the singularity. It  eventually goes away at the boundary $r=R$ of the computational domain.  
We observe the convergence of the solution toward the steady solution using the well balanced second--order Lax--Friedrichs scheme. One observes that the numerical solution, using the well balanced second--order Lax--Friedrichs scheme, is stable. 
Next, in Figure~5.3.2, we treated the second model (\ref{RB5}), with $\Delta r =0.001$ and $v(R)=0.1$, and similarly to the first case, one observes
that the initial perturbation evolves toward a so--called $N$--wave, and that this perturbation moves in the right direction away from the singularity. Again, the solution converges toward the steady solution.

\begin{figure}[t!]   \ 
\centering \subfloat[]{
    \begin{minipage}[b]{0.48\linewidth}
        \centering \includegraphics[height=120pt, width=0.98\linewidth]{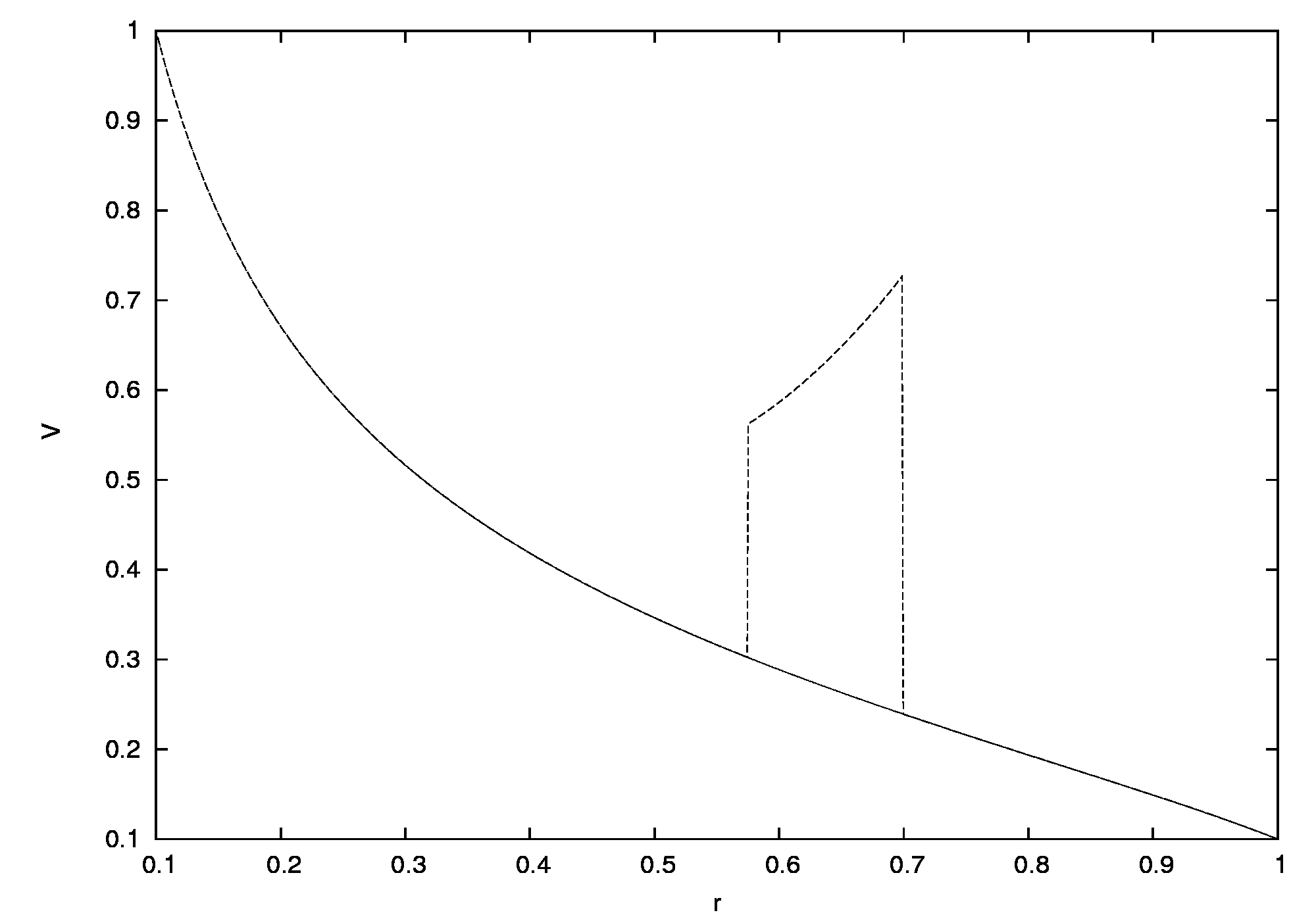}
    \end{minipage}}
\subfloat[]{
    \begin{minipage}[b]{0.48\linewidth}
        \centering \includegraphics[height=120pt, width=0.98\linewidth]{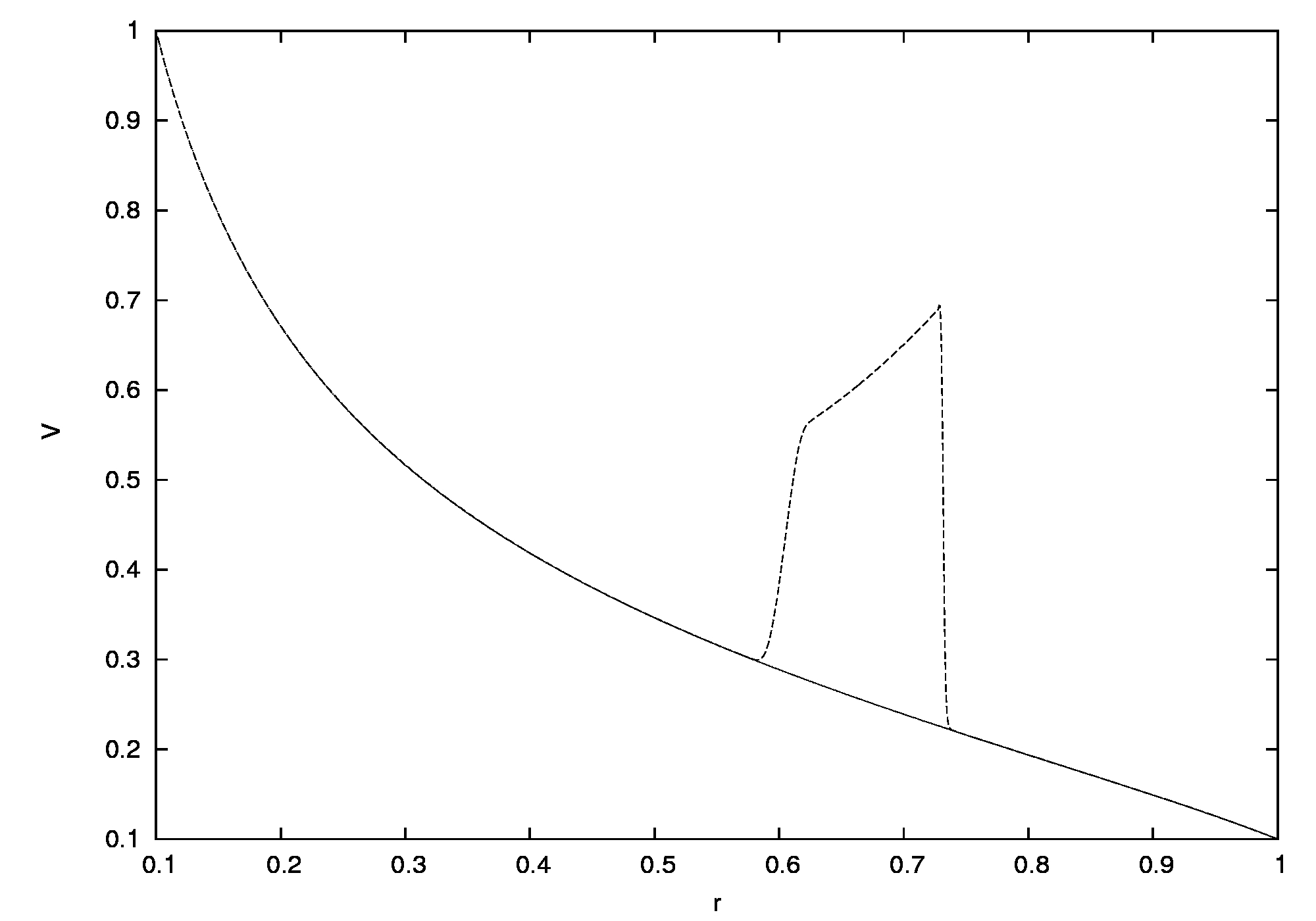}
    \end{minipage}
} \hfill \subfloat[]{
    \begin{minipage}[b]{0.48\linewidth}
        \centering \includegraphics[height=120pt, width=0.98\linewidth]{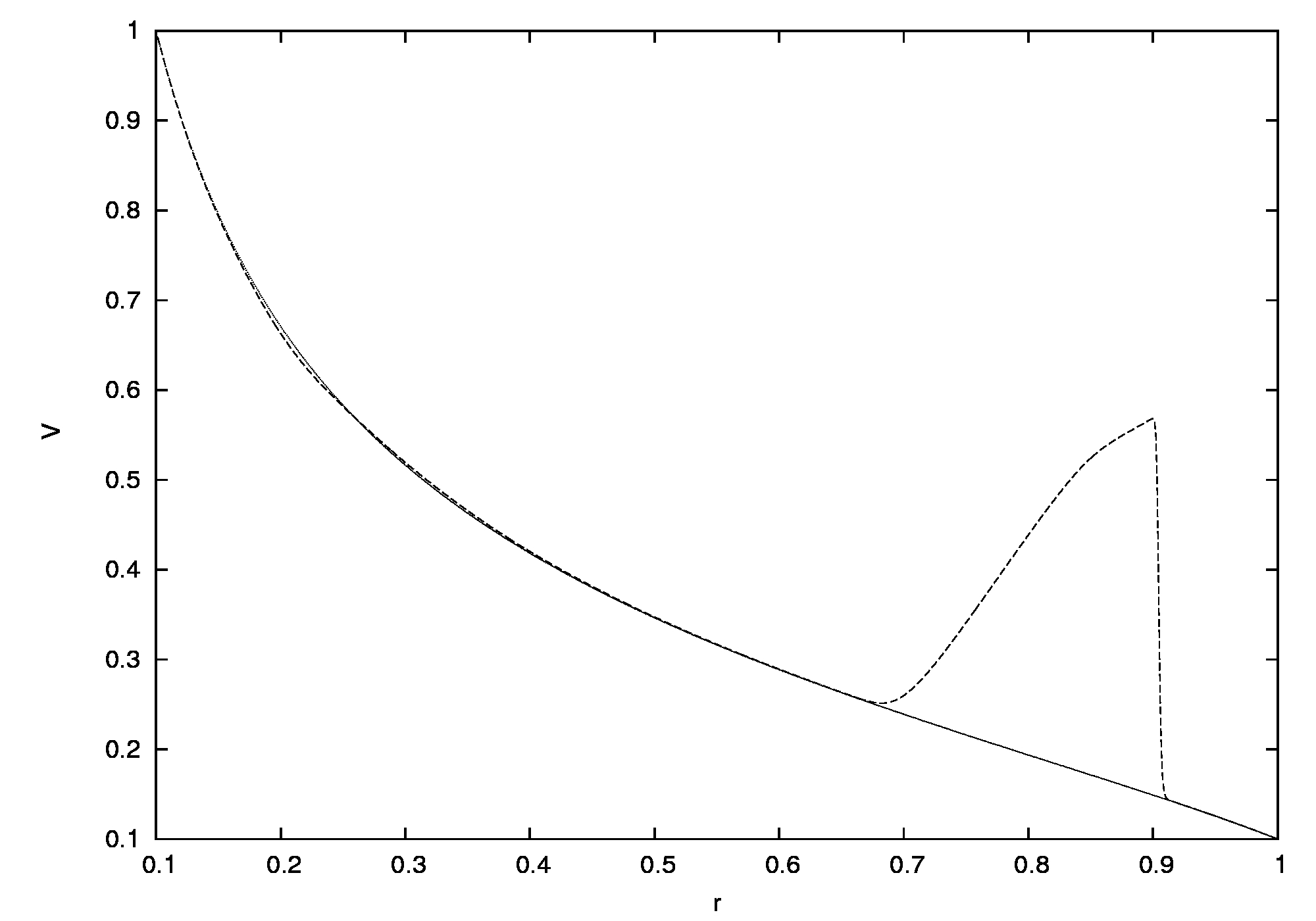}
    \end{minipage}
}
\subfloat[]{
    \begin{minipage}[b]{0.48\linewidth}
        \centering \includegraphics[height=120pt, width=0.98\linewidth]{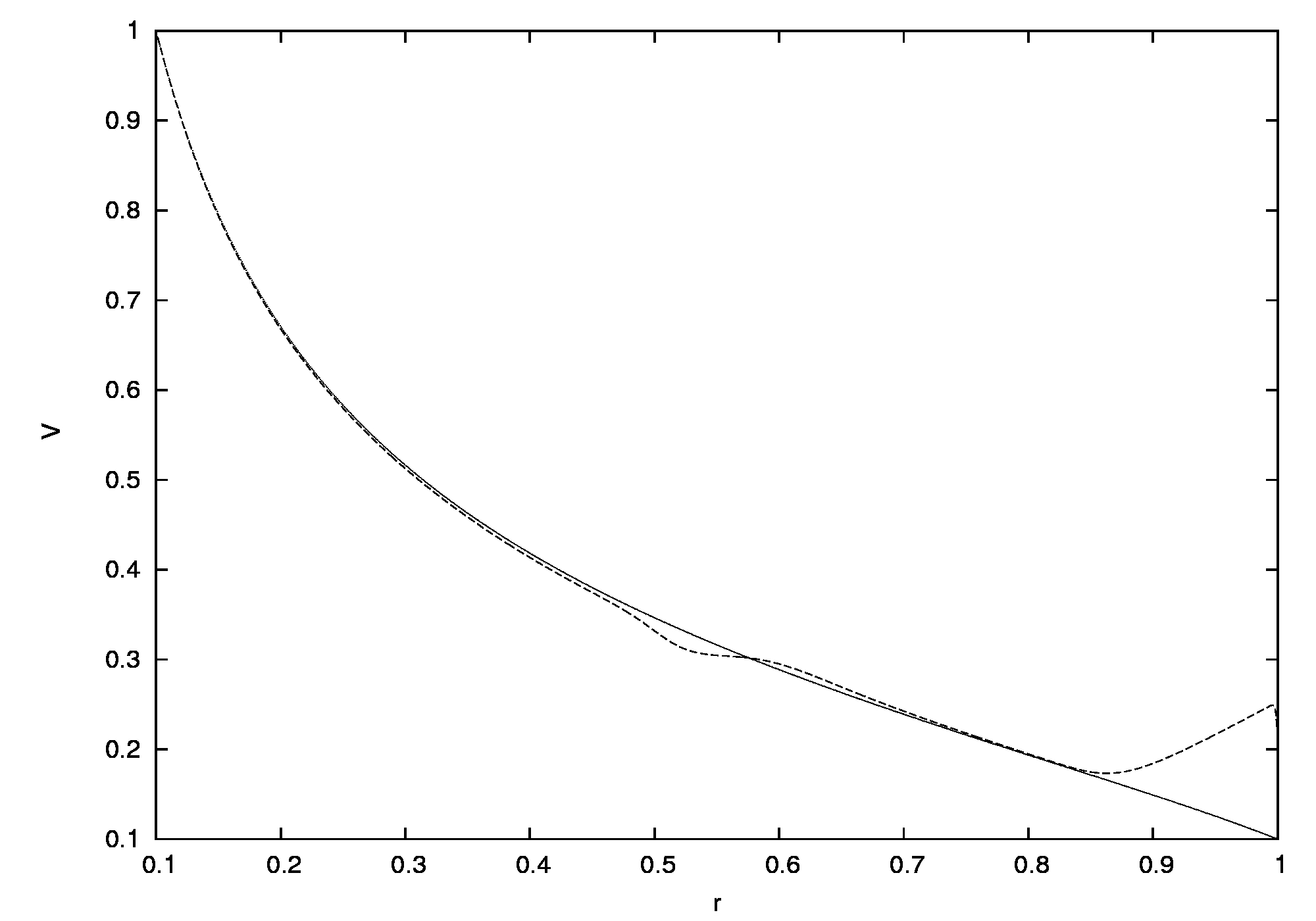}
    \end{minipage}
}\begin{center}
                    {Figure 5.3.2. Model (II). Late-time asymptotic behavior for the model (\ref{RB5}) using the well balanced second--order scheme.}
                   \end{center}
\end{figure}


\subsection{Application to late-time asymptotics -- perturbed one shock solution}

Finally, we take the case of a solution containing a single steady shock separating two initially steady solutions for the equation (\ref{static5}) at the location $r =0.5$. We consider here the first model (\ref{ss5}) with $\Delta r=0.05$ and $v(R)=0.09$.
We study here the nonlinear stability of this steady solution and include 
initial perturbation on a steady solution.
 The numerical results are plotted in Figure~5.4. We observe the convergence of the solution toward the steady solution (single steady shock) using the well balanced second--order Lax--Friedrichs scheme.


\section{Conclusions}

In the first part of this paper, we considered a general class of hyperbolic balance laws posed on a curved spacetimes and we derived flux vector fields $T=T(v)$ which,  in a relativistic context, generalize the classical Burgers flux. 
On one hand, we required that (\ref{AAA}) satisfies the {\sl Lorentz invariance property} shared by the Euler equations of relativistic compressible fluids, and we identified a balance law, which is unique up to normalization. Interestingly enough, the standard Burgers equation was recovered by taking the singular limit of infinite light speed. On the other hand, an alternative model of interest was derived {\sl directly from the relativistic Euler equations,} by assuming that the pressure  term vanishes. This latter model, in the case of a flat background geometry, turned out to be simply the classical Burgers equation. Both models are referred to {\sl relativistic Burgers equations on curved spacetimes.} They exhibit many of the mathematical properties (hyperbolicity, genuine nonlinearity, steady solutions) and challenges (shock wave, late-time asymptotics) encountered with the Euler system of relativistic compressible fluids.
In the second part of this paper, we investigated the numerical discretization of entropy solutions to the above two models. We introduced a finite volume scheme for (\ref{AAA}) defined on curved spacetimes, especially on a background Schwarzschild spacetime. The proposed scheme is fully consistent with the divergence form of the equations and, therefore, applies to weak solutions containing shock waves. More importantly, our scheme is well balanced in the sense that it preserves steady solutions, and extensive numerical experiments demonstrated the convergence of the proposed scheme, and its relevance for computing entropy solutions on a curved background. The generalization to the Euler equations posed on Schwarzschild spacetime are presented in the follow-up work \cite{LeFlochMakhlof}.
 

\end{document}